%% file: 2024_01_18.tex
\documentclass[12pt,oneside]{article}
\newcommand{\Path}{}
\newcommand{\figs}{}
\usepackage{ \Path Paper_style}
\usepackage{ \Path Keywords}
\usepackage{ \Path Category}
\usepackage{ \Path Environments}
\usepackage{tikz-cd}
\usepackage{tabularx} \newcolumntype{L}{>{\raggedright\arraybackslash}X}

\input{Paper_keywords.tex}
\title{The categorical basis of dynamical entropy}

\begin{document}
\author{Suddhasattwa Das\footnotemark[1]}
\footnotetext[1]{Department of Mathematics and Statistics, Texas Tech University, Texas, USA, \url{iamsuddhasattwa@gmail.com} }
\date{\today}
\maketitle
\begin{abstract} Many branches of theoretical and applied mathematics require a quantifiable notion of complexity. One such circumstance is a topological dynamical system - which involves a continuous self-map on a metric space. There are many notions of complexity one can assign to the repeated iterations of the map. One of the foundational discoveries of dynamical systems theory is that these have a common limit, known as the topological entropy of the system.  We present a category-theoretic view of topological dynamical entropy, which reveals that the common limit is a consequence of the structural assumptions on these notions. One of the key tools developed is that of a qualifying pair of functors, which ensure a limit preserving property in a manner similar to the sandwiching theorem from Real Analysis. It is shown that the diameter and Lebesgue number of open covers of a compact space, form a qualifying pair of functors. The various notions of complexity are expressed as functors, and natural transformations between these functors lead to their joint convergence to the common limit. 
\end{abstract}

\section{Introduction} \label{sec:intro}

A bounded dynamical system on a metric space is described by the following simple set of assumptions :

\textbf{Ground assumption.} \emph{ There is a compact metric space $\Omega$ of diameter $1$, and a continuous function $\map : \Omega\to \Omega$. }

Here, we have assumed without any loss of generality that the metric on bounded space $\Omega$ has been normalized to have diameter $1$. The focus of topological dynamical systems theory is to derive and describe properties of the system from these basic assumptions. The objects that are usually in consideration are various invariant behavior such as attractors, invariant sets and omega-limit sets, and various asymptotic properties such as invariant measures, stable and unstable manifolds, and entropy. Our focus is on a notion called entropy. More precisely, we are going to reprove a classical result on entropy using a categorical route. This would involve a categorical reformulation of several topological notions such as open covers, diameter and Lebesgue number. The categorical formulation brings to light various structural properties of the dynamics, i.e. various functors it induces.

\begin{figure}[!ht]\center	\includegraphics[width=.95\linewidth]{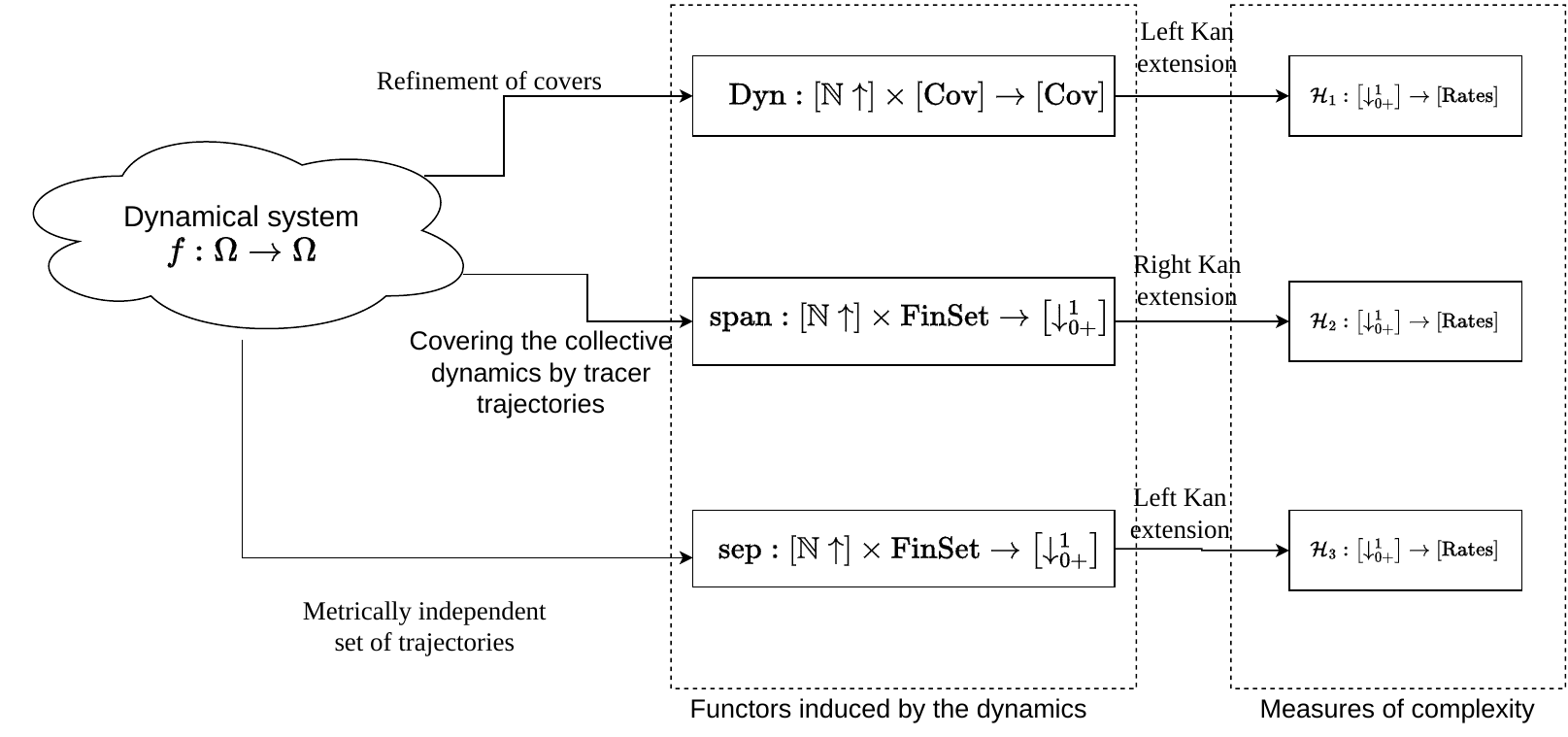}
	\caption{ Outline of the categorical approach. Section~\ref{sec:intro} describes the three classical ways of quantifying the growth of complexity of a dynamical system. We interpret each of these notions of complexity as functors $\dyn$ \eqref{eqn:def:dyn}, $\OrbitSep$ \eqref{eqn:def:OrbCovSep} and  $\OrbitSpan$ \eqref{eqn:def:OrbCovSep}. We have labelled these as the dynamics induced functors. We define a category of growth rates denoted $\Z$ \eqref{eqn:def:rates} and show that the three dynamics induced functors lead to entropy functors $\Pressure_1$ \eqref{eqn:def:P1_P4_fnctr} and $\Pressure_2, \Pressure_3$ \eqref{eqn:def:_P2_P3_fnctr}. These are the functorial representation of the sequences described in Equation \eqref{eqn:pressure}, which describe how the classical notions of complexity grow with successive iterations of the dynamical system.  The asymptotic exponential growth rate of these sequences are interpretable as the colimit functor, and we derive a general categorical result (Theorem \ref{thm:entangle}) to show that they are equal.}
	\label{fig:outline1}
\end{figure}

The role of the dynamics map $f$ is to redistribute the points in $\Omega$. Thus neighborhoods get scattered, entangled or mixed, and the metric also gets finer and less smooth. All these phenomenon reveal in their own way, the complexity of the dynamics.
Topological entropy or simply, entropy, is a measure of the rate of growth of complexity of a dynamical system. There are multiple notions of complexity one can ascribe to a dynamical system, and as a result, multiple definitions of entropy \citep[e.g.][]{Walters1975, Walters2000, PesinPitskel1984, Mummert2007, Climenhaga_entropy, ClimenhagaThompson2012}. They are three main notions of entropy based on partitions, and changes in the metric. All of them are directly inspired from physical phenomena and experimental attempts to observe dynamics. See Figure~\ref{fig:outline1} for a summary of these notions, and our how they are interpreted as functors.

\paragraph{Refinement of partitions} One can track the growth of complexity by tracking how open covers of $\Omega$ get increasingly finer. This notion of complexity is related to the information theoretic complexity associated with data / measurements taken from the dynamics of $\paran{ \Omega, f }$ \citep[e.g.][]{grassberger1991inf, vallee1987inf}, and has important consequences to computability \cite{Spandl_sofic_2008, BDWY2020}. Suppose one intends to track the dynamics on the phase space $\Omega$ by means of a computer. Being a finite state machine, the computer approximates $\Omega$ via a discrete collection of sets whose union covers $\Omega$. To avoid boundary issues, we assume that these sets are open sets. Thus, the matter of observation of a space $\Omega$ becomes naturally associated with the concept of an \emph{open cover}.

Given an open cover $\alpha$ of $\Omega$, any point $x$ on $\Omega$ receives an address-identity based on whichever pieces of $\alpha$ it lies in. If $x$ lies at the intersection of multiple pieces, this address is obviously non-unique. Moreover, different points in $\Omega$ may be assigned the same address-identity.  Recall that the \emph{diameter} of a set is the maximum possible distance between any two points within it. If all the pieces of $\alpha$ have a diameter less than $\epsilon$, then any address-identity of $x$ specifies its location to within an error bound of $\epsilon$.

Dynamical systems is about the study of trajectories. Continuing the same observation scheme, the $n$-length orbit of $x$, which is the sequence $x, f(x), \ldots, f^{n-1}(x)$, is specified by a sequence of address locations. In other words, an orbit is to be described by the computer as a sequence $\alpha_{i_0} , \ldots, \alpha_{i_{n-1}}$ of pieces of $\alpha$, such that $f^t(x) \in \alpha_{i_t}$. Now note that given this sequence $\alpha_{i_0} , \ldots, \alpha_{i_{n-1}}$, the location of the starting point $x$ gets narrowed down to the set 
\[ \alpha \paran{ i_0, \ldots, i_{n-1} } := \cap_{t=0}^{n-1} f^{-t} \alpha_{i_t} . \]
Thus with each iteration of the dynamical system, and each additional effort to locate the new position of the point $x$, one gets a progressively refined estimate of the initial location of $x$. In fact the set $\alpha \paran{ i_0, \ldots, i_{n-1} }$ above is a piece of the \emph{dynamically refined} partition
\[ \alpha^{(n)} = \vee_{t=0}^{n-1} f^{-t} \alpha = \SetDef{ \alpha \paran{ i_0, \ldots, i_{n-1} } }{ \mbox{all sequences } i_0, \ldots, i_{n-1} } . \]
But now note that for an arbitrary sequence $i_0, \ldots, i_{n-1}$, the intersection
$\alpha \paran{ i_0, \ldots, i_{n-1} }$ may be empty. Moreover, one may not need all of the pieces of $\alpha^{(n)}$ to cover $\Omega$. From a computational point of view, one would need to retain a minimum number $\CovNum \paran{ \alpha^{(n)} }$ of pieces. This number $\CovNum \paran{ \alpha^{(n)} }$ may be interpreted as the minimum size of a subcover of $\alpha^{(n)}$. Note that $\CovNum \paran{ \alpha^{(n)} }$ is finite by the compactness of $\Omega$.

In summary, any open cover $\alpha$ (finite or not) produces a sequence of numbers
\[ \CovNum \paran{ \alpha^{(1)} }, \CovNum \paran{ \alpha^{(2)} }, \ldots , \CovNum \paran{ \alpha^{(n)} }, \ldots \]
The first notion of entropy studies the asymptotic exponential rate at which this sequence grows, as the size / diameter of the partition $\alpha$ diminishes. The open cover $\alpha$ can be interpreted to be a coarse grained partition of division of the space $\Omega$. Then  $\CovNum \paran{ \alpha^{(n)} }$ is the minimum storage required to use the alphabet / addresses of $\alpha$ to describe all possible trajectories of length $n$.

\paragraph{Refinement of the metric} Instead of tracking the dynamics using a cover, one could release a finite collection of tracers or sensors into the phase space, and track their movements. Such tracers can be specified by their set of initial locations, a finite subset $A$ of $\Omega$. Each $a\in A$ has a trajectory $\{f^t(a)\}_{t=0,1,2,\ldots}$, and the objective is that an arbitrary trajectory $\{f^t(x)\}_{t=0,1,2,\ldots}$ can be followed or closely approximated by at least one of the tracer trajectories. More precisely, for each $n$ and finite set $A$, we are interested in
\begin{equation} \label{eqn:odp92}
    \sup_{x\in\Omega} \; \inf_{a\in A} \; \max_{0\leq t<n} \; d \paran{ f^t x, f^t a } .
\end{equation}
The dual to this problem of covering the phase spaces with tracers, is the task of finding independent trajectories. The independence of the points in a finite set $A$ can be measured as the closest their trajectories get to each other, i.e., as the quantity
\begin{equation} \label{eqn:pdmx34}
    \inf_{a\neq a'\in A} \; \max_{0\leq t < n} d \paran{ f^t a, f^t a' } .
\end{equation}
The choice of such a set of independent points is crucial in applications such as finding \emph{landmark points} \citep[e.g.][]{LongFerguson2019, LiangPaisley2015, SilvaEtAl_landmark_2006} for manifold learning. The min-max problems of \eqref{eqn:odp92} and \eqref{eqn:pdmx34} may be understood better from the realization that the starting metric $d$ on $\Omega$ is refined under repeated action of the dynamics. For each $n$, define
\[ d_n (x,y) = \max_{0\leq t<n} d \paran{ f^t x, f^t y } .\]
This leads to  a series of metrics $d = d_0 << d_1 << d_2 << \ldots$, where  $<<$ indicates that each metric is successively finer. The metric $d_n$ incorporates some information about the dynamics, by taking into account the future $n$ states of each point. Given any metric $\tau$ on $x$ that generates its topology, and finite set $A\subset \Omega$ we can define two quantities
\begin{equation} \label{eqn:dpd9}
    \text{metric-sep} (A,\tau) := \min \SetDef{ \tau(x,y) }{ x,y \in A, \, x\neq y } . \quad \text{metric-span} (A,\tau) := \max \SetDef{ \tau(x,A) }{ x\in X } .
\end{equation}
As $n$ grows, the metric $d_n$ grows finer. As a result, both the quantities $\text{metric-sep}$ and $\text{metric-span}$ increase.  As a result, $\text{metric-span}$ and $\text{metric-sep}$ can reformulate the optimization problems in \eqref{eqn:odp92} and \eqref{eqn:pdmx34} as
\begin{equation} \label{eqn:def:P2_P3} 
    \begin{split}
        H_2( \epsilon, n) &:= \min \SetDef{ \abs{A} }{ A \subset_{fin} \Omega, \; \text{metric-span} \paran{ A, d_n } > \epsilon } . \\
        H_3( \epsilon, n) &:= \max \SetDef{ \abs{A} }{ A \subset_{fin} \Omega, \; \text{metric-sep} \paran{ A, d_n } < \epsilon } .
    \end{split} 
\end{equation}
Thus, $H_2( \epsilon, n)$ is the smallest size of a discrete set $A$ such that the $n$-orbit of every point is within $\epsilon$ distance of the $n$-orbit of some point in $A$. Similarly, $H_3$ covers a dual notion of how separated $n$-orbits remain. 

\paragraph{Goal} The following well known result from ergodic theory \citep[e.g.][Theorem 9.4]{Walters2000} shows that complexity may be measured in $4$ different ways, and each have the same asymptotic limit, called the \emph{topological entropy} of the dynamics :
\begin{corollary} [Topological entropy] \label{cor:thermo}
For any continuous dynamical system on a bounded metric space, 
\begin{equation} \label{eqn:pressure}\begin{split}
& \lim_{\epsilon\to 0^+} \sup_{ \diam(\alpha) \geq \epsilon } \lim_{n\to\infty} \frac{1}{n} \ln \CovNum \paran{ \alpha^{(n)} } = \sup_{ \alpha } \lim_{n\to\infty} \frac{1}{n} \ln \CovNum \paran{ \alpha^{(n)} } \\
& = \lim_{\epsilon\to 0^+} \lim_{n\to\infty} \frac{1}{n} \ln H_2(\epsilon, n) = \lim_{\epsilon\to 0^+} \lim_{n\to\infty} \frac{1}{n} \ln H_3( \epsilon, n ) .
\end{split}\end{equation}
\end{corollary}
Our goal is to rediscover Corollary~\ref{cor:thermo} by interpreting several notions of dynamics, entropy and limits using categorical means. All the limits in \eqref{eqn:pressure} will be interpreted as colimits of functors, as summarized in Figure~\ref{fig:outline1}. Their equality will be established by revealing a web of relations between the functors, and is based on a purely categorical result we state and prove as Theorem \ref{thm:entangle}. An overview of these inter-relations is also presented in Figure~\ref{fig:outline2}. These relations are based on an interpretation of some basic topological constructs, namely covers, diameter, and Lebesgue number. We show that these familiar notions can be reinterpreted as functors. We derive two other purely category theoretical results in Theorems \ref{thm:qualify} and \ref{thm:lim_qual}, which brings to light the structural role played by these concepts. Corollary~\ref{cor:thermo} eventually turns out to be a consequence of these reinterpretations of topology and dynamics.

\paragraph{Previous work.} Due to the many facets of dynamical systems, there are different ways one can find categorical structure in dynamics. For example, one could characterize the dynamics via its limiting sets \cite{calcines2013limit}, their state-space action \cite{spivak2015steady, ngotiaoco2017}, collection of orbits \cite{jaz2020double}. There has been notable investigations into  entropy in the context of measure theory \citep[e.g.][]{Baez2014_bayesian, BaezEtAl2011_info} separately. Most of the categorical treatment of dynamical systems \citep[e.g.][]{lomadze1999time, Delvenne2019_dyn, Suda2022Poincare, MossPerrone2022ergdc} deal with the global structure of dynamical systems and do not pursue the inner, topological changes that take place. Our work provides a new joint perspective into dynamics and topology. An important consideration for us the existence and equalities of limits, which from a categorical point of view, translates into colimits. One of the key contributions is the notion of a \emph{qualifying pair} of functors, defined and studied in Section~\ref{sec:qualify}. We show how such a pair can be created simultaneously, and have important consequences to the existence of colimits. More importantly, they play a central role in the existence and uniqueness of the limits in Corollary~\ref{cor:thermo}.

\paragraph{Outline} We begin our categorical approach by re-examining some basic concepts of topology in Section~\ref{sec:topo}. The most important outcome is that the notions of diameter and Lebesgue number of a partition are reinterpreted as functors, or more precisely, as \emph{Kan extensions}. In Section~\ref{sec:qualify} we introduce the notion of a qualifying pair. Diameter and Lebesgue number are shown to be such a pair. The existence of qualifying pairs have important consequences when studying limiting behavior, or more precisely, colimits. We examine this in detail in Section~\ref{sec:qpl}. This completes the category theoretical background for the work. Next, we look at the functors that are induced by the dynamics in Section~\ref{sec:dynamics}. The notions of complexity introduced in \eqref{eqn:def:P2_P3} and \eqref{eqn:pressure} are given a functorial redefinition. Finally, Corollary~\ref{cor:thermo} is proved in Section~\ref{sec:network}. This is done by interlinking the complexity functors using the category theoretical tools we have developed. The proofs of the lemmas and Theorems that are used on the way, are provided in Section~\ref{sec:proofs}.

\paragraph{Notation} Throughout the rest of the paper, we will use the notation $F:X\to Y$ to denote a functor between categories $X, Y$. Given two functors $F,G : X\to Y$, the notation $F\Rightarrow G$ signifies a natural transformation from $F$ to $G$. Identity map. The set of objects of a category $\mathcal{C}$ will be denoted as $ob\paran{\mathcal{C}}$. 

\section{Building blocks of metric topology} \label{sec:topo}

One of the fundamental notions of topology is that of open covers. Open covers are a means of granularizing or discretizing a continuum. Any finite representation of a continuum, such as using a computer, relies on the creation of an open cover. The complexity of a dynamical system may be studied by means of such representation using open covers. We shall now discover that the set of open covers $\OpenCov$ form a preorder, i.e., a category in which there is at most one morphism between any two objects. 

Given a preorder, a morphism $f:a\to b$ between two objects can be unambiguously denoted by an arrow $a\to b$. Thus the morphisms together create a partial ordering of the objects. There are two ways one can assign a partial ordering to $\OpenCov$. For any two open covers $a, b$ :
\begin{enumerate}[(i)]
    \item $a\to b$ if for every set $S_b$ in $b$, there is a set $S_a$ in $a$ such that $S_a \supseteq S_b$.
    \item $a\to b$ if for every set $S_a$ in $a$, there is a set $S_b$ in $b$ such that $S_a \supseteq S_b$.
\end{enumerate}
These preorder structures on the collection of open covers are different and provide different interpretations of the notion - `` $a$ is coarser than $b$" via the morphism $a\to b$. We choose the first option. We next consider a more trivial subclass of open covers : the collection $\udc$ of covers all of whose pieces are open disks of a constant radius. We call this the subcategory of \emph{uniform disk covers}. 

There are a few other fundamental preorder categories which are tied to the study of $\OpenCov$. The first is $\Rplus$ - the preorder of numbers in $[0,\infty]$ ordered by the $\leq$ relation. Its opposite category is $\Rminus$ - the preorder of numbers in $[0,\infty]$ ordered by the $\geq$ relation. One can similarly define $\Nplus$ and $\Nminus$, the subcategories of $\Rplus$ and $\Rminus$ respectively, generated by the integers. Finally, we define $\EpsCat$ to be the subcategory of $\Rminus$ obtained by restriction to the set $(0,1]$. It is useful to model any process of \emph{parameterization} by a bounded, non-negative number $\epsilon$ as $\EpsCat$, if one is interested in the limiting behaviour as $\epsilon\to 0^+$. 

The categories $\udc$, $\OpenCov$ and $\EpsCat$ are linked together by some functors fundamental to metric topology. First, there is an inclusion map $\iota_{ UDC } : \udc \to \OpenCov$ as an inclusion of subcategories. Next, recall that a functor between two preorders is the same as an order preserving map between the respective set of objects. As a result, the map which assigns every uniform disk cover in $\udc$ its constant radius, is also a functor 
\[ \Grain : \udc \to \EpsCat \]
The $\Grain$ functor might seem trivial, but it is the key to building the more universal idea of diameter. Next, one can assign to each $\epsilon \in ob\paran{ \OpenCov }$ the open cover containing all possible disks of radius $\epsilon$. This is the free UDC generated by $\epsilon$. One in fact can represent this as a functor 
\[ \Free : \EpsCat \to \udc . \]
The categories $\udc, \OpenCov, \EpsCat$ and the functors $\Grain, \Free$ will be the building blocks of topological entropy. The study of open covers and open sets involve the use of several quantifiable notions such as diameter and Lebesgue number. The advantage of the categorical point of view is that these notions can be derived as structural properties. We next describe \emph{Kan extensions}, one of the key tools to derive new structural properties from an existing arrangement of categories and functors. 

\paragraph{Kan extensions} Kan extensions \citep[e.g.][]{perrone2022kan, street2004categorical, Riehl_homotopy_2014} are universal constructions which generalize the practice of taking partial minima or maxima, in a functorial manner. Suppose we have two functors $E \xleftarrow{F} X \xrightarrow{K} D$. Although not necessary to the definition, one can interpret $F$ as a measurement on a parameter space $X$, and $K:X\to D$ to be a projection to a subset of coordinates. A \emph{left Kan extension} of $F$ along $K$ is a functor $\REnv{K}{F} : D\to E$ along with a minimum natural transformation $\eta : F\Rightarrow \REnv{K}{F} \circ K$. Moreover, this pair $\paran{ \REnv{K}{F}, \eta }$  is minimum / universal in the sense that for every other functor $H:D\to E$ along with a natural transformation $\gamma : F \Rightarrow H\circ K$, there is a natural transformation $\tilde\gamma : \REnv{K}{F} \Rightarrow H$ s.t. $\gamma = \paran{ \tilde\gamma \star \Id_{K} } \circ \eta$.  This is shown in the diagram below.
\[ \mathcal{L} := \REnv{K}{F}, \quad 
\begin{tikzcd}
    & & & E \\
    E & X \arrow[dashed, bend left = 10]{urr}[name=x1]{} \arrow[dashed]{drr}[name=x2]{} \arrow{l}[name=F]{F} \arrow{r}{K} & D \arrow{dr}[name=H]{H} \arrow{ur}[swap, name=L]{ \mathcal{L} } & \\
    & & & E 
    \arrow[shorten <=2pt, shorten >=3pt, Rightarrow, to path={(F) to[out=90,in=180] (x1)} ]{  }
    \arrow[shorten <=2pt, shorten >=3pt, Rightarrow, to path={(F) to[out=-90,in=225] (x2)} ]{  }
    \arrow[shorten <=1pt, shorten >=1pt, Rightarrow, to path={(L) to[out=-45,in=45] (H)} ]{  }
\end{tikzcd}\]
One can similarly define a \emph{right Kan extension} of $F$ along $K$. It is a functor $\LEnv{K}{F} : D\to E$ along with a natural transformation $\epsilon : \LEnv{K}{F} \circ K \Rightarrow F$. Moreover, this pair $\paran{\LEnv{K}{F}, \epsilon}$ is maximum / universal in the sense that for every other functor $H:D\to E$ along with a natural transformation $\gamma : H\circ K \Rightarrow F$, there is a natural transformation $\tilde{\gamma} : H \Rightarrow \LEnv{K}{H}$ such that $\gamma = \epsilon\circ \paran{ \tilde{K} \star \Id_K }$. This is shown in the diagram below.
\[ \mathcal{R} := \LEnv{K}{F}, \quad 
\begin{tikzcd}
    & & & E \\
    E & X \arrow[dashed, bend left = 10]{urr}[name=x1]{} \arrow[dashed]{drr}[name=x2]{} \arrow{l}[name=F]{F} \arrow{r}{K} & D \arrow{dr}[name=H]{H} \arrow{ur}[swap, name=L]{ \mathcal{R} } & \\
    & & & E 
    \arrow[shorten <=2pt, shorten >=3pt, Rightarrow, to path={(x1) to[out=135,in=90] (F)} ]{  }
    \arrow[shorten <=2pt, shorten >=3pt, Rightarrow, to path={(x2) to[out=225,in=-90] (F)} ]{  }
    \arrow[shorten <=1pt, shorten >=1pt, Rightarrow, to path={(H) to[out=45,in=-45] (L)} ]{  }
\end{tikzcd}\]
If $E$ is a co-complete category, the left Kan extension always exists. Similarly if $E$ is a complete category, the right Kan extension always exists. In case both left and right Kan extensions of $F$ along  $K$ exist, they combine to produce the following diagram :
\begin{equation} \label{eqn:L_R_Env}
    \begin{tikzcd} [column sep = large]
        E & E & E \\
        & X \arrow[dashed]{ul}[name=x2]{} \arrow[dashed]{ur}[name=x1]{} \arrow{d}{K} \arrow{u}[name=F]{F} \\ 
        & D \arrow[]{uur}[swap]{\REnv{K}{F}} \arrow[]{uul}{\LEnv{K}{F}} &
        \arrow[shorten <=2pt, shorten >=3pt, Rightarrow, to path={(F) -- (x1)},]{}
        \arrow[shorten <=2pt, shorten >=3pt, Rightarrow, to path={(x2) -- (F)},]{}
    \end{tikzcd}
\end{equation}
We prove in Lemma~\ref{lem:cov} that $\Grain,\Free$ satisfy the following assumptions : 

\begin{Assumption} \label{A:1}
There is a embedding $\iota_{X} : \tilde{X}\to X$ of categories. 
\end{Assumption}

\begin{Assumption} \label{A:3}
There are functors $G : \tilde{X}\to \mathcal{E}$ and $G^r : \mathcal{E}\to \tilde{X}$ such that $G G^{r} = \Id_{\mathcal{E}}$ and $G^{r} G \Rightarrow \Id_X$.
\end{Assumption}

\begin{Assumption} \label{A:2}
The category $\mathcal{E}$ is complete, and the right Kan extension $\LEnv{ \Id_{ \mathcal{E} } }{ \iota \Free }$ of $\Id_{ \mathcal{E} }$ along $\iota G^r$ exists
\end{Assumption}

\begin{lemma} \label{lem:cov}
Assumptions \ref{A:1}, \ref{A:2} and \ref{A:3} are satisfied if one takes $\udc, \OpenCov, \EpsCat, \Grain$ and $\Free$ to be $\tilde{X}, X, \mathcal{E}, G$ and $G^{r}$ respectively.
\end{lemma}

The totally ordered set $\EpsCat$ is clearly complete. Covers in $\udc$ are clearly open covers and thus elements in $\OpenCov$. Thus the inclusion of $\udc$ in $\OpenCov$ is an embedded $\iota_{udc}$. Since all the categories involved are preorders, and the maps $\Grain$ and $\Free$ are order preserving, $\Grain$ and $\Free$ are also functors. To verify Assumption~\ref{A:2}, note that the completeness of $\EpsCat$ stems from the fact that its objects form the right closed interval $(0,1]$. The left Kan extension exists and is defined pointwise as
\[ \LEnv{ \Id_{ \mathcal{E} } }{ \iota \Free }(\alpha) := \sup \SetDef{ \epsilon\in (0,1] }{ \alpha \to \Free(\epsilon} .\]
Recall that the Lebesgue number of an open cover $\alpha$ is the maximum number $\delta$ such that for every $\delta'<\delta$ and any ball $B$ of diameter $\delta'$, there is a piece of $\alpha$ which contains $B$. In fact Lebesgue number is a Kan extension, as stated below :

\begin{lemma} \label{lem:Leb_diam}
The Lebesgue number and diameter of an open cover are right and left Kan extensions :
\[ \LebNum := \LEnv{ \Id_{ \mathcal{E} } }{ \iota \Free } : \OpenCov \to \EpsCat, \quad \diam := \REnv{ \Grain }{ \iota_{UDC} } : \OpenCov \to \EpsCat. \]
\end{lemma}

Lemma~\ref{lem:Leb_diam} is proved in Section~\ref{sec:proof:Leb_diam}. the constructions of $\LebNum$ and $\diam$ can be stated more abstractly using the notations of Assumptions \ref{A:1}, \ref{A:3}. Thus, one can similarly define the left and right Kan extensions
\begin{equation} \label{eqn:def:LD}
    L := \LEnv{ \Id_{ \mathcal{E} } }{ \iota G^r } : X \to \mathcal{E}, \quad D := \REnv{G}{ \iota } : X \to \mathcal{E} .
\end{equation}
The commutation between the various spaces are shown on the left part of \eqref{eqn:doer4}. The right part relates this abstract framework to our concrete examples in metric topology.
\begin{equation} \label{eqn:doer4}
    \begin{tikzcd} 
        & & \mathcal{E} & \\
        \mathcal{E} & \mathcal{E} \arrow{ur}[name=ide]{ \Id_{\mathcal{E}} } \arrow[dashed]{l}[name=L]{} \arrow{r}{ G^r } & \tilde{X} \arrow[dashed]{r}[name=D]{} \arrow{d}{\iota} \arrow{u}[swap, name=G]{G} & \mathcal{E} \\
        & & X \arrow{llu}{ L := \LEnv{ \Id_{ \mathcal{E} } }{ \iota G^r } } \arrow{ru}[swap]{ D := \REnv{G}{\iota} } & 
        \arrow[shorten <=2pt, shorten >=3pt, Rightarrow, to path={(L) to[out=90,in=180] (ide)} ]{}
        \arrow[shorten <=2pt, shorten >=3pt, Rightarrow, to path={(G) to[out=0,in=90] (D)} ]{}
    \end{tikzcd}, \quad
    \begin{tikzcd} 
        & & \EpsCat & \\
        \EpsCat & \EpsCat \arrow{ur}[name=ide]{ \Id } \arrow[dashed]{l}[name=L]{} \arrow{r}{ \Free} & \udc \arrow[dashed]{r}[name=D]{} \arrow{d}{\iota} \arrow{u}[swap, name=G]{ \Grain } & \EpsCat \\
        & & \OpenCov \arrow{llu}{ \LebNum } \arrow{ru}[swap]{ \diam } & 
        \arrow[shorten <=2pt, shorten >=3pt, Rightarrow, to path={(L) to[out=90,in=180] (ide) } ]{}
        \arrow[shorten <=2pt, shorten >=3pt, Rightarrow, to path={(G) to[out=0,in=90] (D)} ]{}
    \end{tikzcd}
\end{equation}
So far we have presented the elementary notions of open covers and disk-open covers as categories $\OpenCov$ and $\udc$, and then shown how the more sophisticated notions of diameter and Lebesgue number are consequences of the mere arrangements between $\OpenCov$ and $\udc$. In fact, the technique of Kan extensions leads to a more significant property for the pair $\paran{ \diam, \LebNum }$. We discuss this next.

\section{Qualifying pair of functors} \label{sec:qualify}

The diameter and Lebesgue number functors can be arranged as $\OpenCov \xrightarrow{\LebNum} \EpsCat \xleftarrow{\diam} \OpenCov$. Recall that given a general arrangement of functors and categories such as $D \xrightarrow{F} C \xleftarrow{G} E$, one can create a \emph{comma category} \citep[e.g.][]{polonsky2020local} $\Comma{F}{G}$ whose objects are 
\[ \SetDef{ \paran{ d, e, \phi } }{ d\in ob(D), \, e\in ob(E), \, \phi \in \Hom_C \paran{ Fd; Ge } } . \]
A morphism between two objects $(d,\phi,e)$ and $(d',\phi',e')$ in this category are all those pairs of morphisms $\SetDef{ (f,g) }{ f\in \Hom_D \paran{d;d'},\, g\in \Hom_E \paran{ e;e' } }$ characterized by the following equation :
\[ 
\begin{tikzcd} (d,\phi,e) \arrow{r}{ (f,g) } & (d',\phi',e') \end{tikzcd} \,\Leftrightarrow \,
    \begin{tikzcd} d \arrow{d}{f} \\ d' \end{tikzcd},  \begin{tikzcd} e \arrow{d}{g} \\ e' \end{tikzcd}, \mbox{ and }
\begin{tikzcd}
    Fd \arrow{d}{Ff} \arrow{r}{\phi} & Ge \arrow{d}{Gg} \\
    Fd' \arrow{r}{\phi'} & Ge'
\end{tikzcd}
\]
Comma categories thus store the commutation relations between $F$ and $G$, when projected to $C$. Comma categories have a \emph{forgetful functor}
\[ \Forget : \Comma{F}{G} \to D\times E , \]
which forgets the commutation information stored within the comma category and projects onto the base category. The action of $\Forget$ on objects and morphisms is described as
\[ (d,\phi,e) \stackrel{ \Forget }{ \mapsto } (d,e), \quad 
\begin{tikzcd} (d,\phi,e) \arrow{d}{ (f,g) } \\ (d',\phi',e') \end{tikzcd} \stackrel{ \Forget }{ \mapsto } 
\begin{tikzcd} (d,e) \arrow{d}{ f\times g } \\ (d',e') \end{tikzcd}
\]
In the special case when $D=E = C$ and $F=G=\Id_C$, the resulting comma category $\Comma{\Id_C}{\Id_C}$ is known as the \emph{arrow category} \citep[e.g.][]{winter2009arrow, steingartner2014categorical} of $C$. It is the category arising out of the commutation relations in $C$. Comma categories arise naturally in many applications of category theory 
\citep[e.g.][]{grandis1997categorically, burke2018synthetic}, and even in the pointwise construction of Kan extensions. We use the language of comma and arrow categories to define a special property for a pair of functors within the same functor category.

\paragraph{Qualifying pairs} Let $C,D$ be categories. A pair of functors $F,G:C\to D$ is said to be a \emph{qualifying pair} if $F\Rightarrow G$, and there is a functor $\Psi : \Comma{F}{G} \to \Comma{ \Id_{C} }{ \Id_{C} }$ satisfying the following commutation of functors :
\[\begin{tikzcd}
    \Comma{F}{G} \arrow[bend left = 20]{rr}{ \Psi } \arrow{r}[swap]{ \Forget } & C \times C & \Comma{ \Id_{ C } }{ \Id_{ C } } \arrow{l}{ \Forget }
\end{tikzcd}\]
An object $\paran{c_1, c_2, \phi }$ of $\Comma{F}{G}$ indicates the morphism $Fc_1 \xrightarrow{\phi} Gc_2$ in $D$. This object is mapped into $\paran{c_1, c_2}$ under $\Forget$. The presence of the commuting map $\Psi$ which is equivariant under $\Forget$ means that $Fc_1 \xrightarrow{\phi} Gc_2$ is mapped into $c_1 \xrightarrow{ \Psi(\phi) } c_2$. Moreover, this mapping is functorial, meaning that morphisms are mapped into morphisms. This is illustrated in the diagram below, in which the red and blue arrows are corresponding objects in their respective comma categories.
\[\begin{tikzcd}
    \red{Fc_1} \arrow[red]{r}{\phi} \arrow{d}[swap]{Ff_1} & \red{Gc_2} \arrow{d}[name=x1]{Gf_2} & & & \red{ c_1 } \arrow[red]{r}{ \Psi(\phi) } \arrow{d}[swap, name=x2]{f_1} & \red{c_2} \arrow{d}{f_2} \\
    \blue{Fc'_1} \arrow[blue]{r}[swap, name=x4]{\phi'} & \blue{Gc'_2} & & & \blue{c'_1} \arrow[blue]{r}[swap, name=x5]{ \Psi(\phi') } & \blue{c'_2} \\
    & & & \red{ \paran{c_1 , c_2} } \arrow{d}[name=x3]{ f_1, f_2 } \\
    & & & \blue{\paran{c'_1 , c'_2}}
    \arrow[shorten <=4pt, shorten >=4pt, Rightarrow, to path={(x1) -- (x2)},]{}
    \arrow[shorten <=4pt, shorten >=4pt, Rightarrow, to path={(x4) -- (x3)},]{}
    \arrow[shorten <=4pt, shorten >=4pt, Rightarrow, to path={(x5) -- (x3)},]{}
\end{tikzcd} \]
The first main result shows that the basic Assumptions  \ref{A:1}, \ref{A:3} guarantee the existence of a qualifying pair.

\begin{theorem}[Qualifying pair construction] \label{thm:qualify}
Let Assumptions \ref{A:1}, \ref{A:3} hold, and assume the notation of \eqref{eqn:def:LD}. Then 
\begin{enumerate}[(i)]
    \item $D$ and $L$ are left and right Kan extensions respectively, of $\Id_{\mathcal{E}}$ along $\iota G^r$.
    \item $D \Rightarrow L$.
    \item Now suppose that $\mathcal{E}$  is totally ordered. Then $(L,D)$ form a qualifying pair of functors.
\end{enumerate}
\end{theorem}

Theorem~\ref{thm:qualify} is proved in Section~\ref{sec:proof:qualify}. Theorem~\ref{thm:qualify}~(i) thus offers a definition for $D$ and $L$ which is an alternative to \eqref{eqn:def:LD}. Now they are revealed to be left and right Kan extensions respectively, of $\Id_{\mathcal{E}}$ along $\iota G^r$. By virtue of Lemma~\ref{lem:cov}, the fundamental metric notions of diameter and Lebesgue number now appear to be left and right Kan extensions respectively along a trivial inclusion functor. 

The next result reveals the rationale behind the name ``qualifying pair". 

\begin{theorem}[Qualifying action] \label{thm:qlf_actn}
    Let $I,C,D$ be categories and $F,G : C\to D$ be a qualifying pair of functors. Then for any pair of functors $F':I\to C$ and $G':I\to C$, any natural transformation $\eta : FF' \Rightarrow GG'$ is mapped into a natural transformation $\Psi_* \eta : F'\Rightarrow G'$.
\end{theorem}

Theorem~\ref{thm:qlf_actn} is proved in Section~\ref{sec:proof:qlf_actn}. Thus for any pair of functors $F',G'$ from the functor category $\Functor{I}{C}$, the composite $FF'$ can be naturally transformed to $GG'$ only if $F'$ can be naturally transformed to $G'$. By virtue of Lemma~\ref{lem:cov} and Theorem~\ref{thm:qualify}, this can be interpreted in the context of metric  topology as follows :

\begin{corollary} [Diameter-Lebesgue num. pair] \label{cor:sdol} 
    The pair of functors $\paran{ \LebNum, \diam }$ is a qualifying pair. Moreover, for any category $I$ and functors $F',G':I \to \OpenCov$, 
    \[ \LebNum( F'(i) ) \leq \diam( G'(i) ), \; \forall i\in ob(I) \quad \Leftrightarrow \quad F'(i) < G'(i), \; \forall i\in ob(I) . \]
\end{corollary}

Thus Lebesgue number of a cover $\alpha$ being smaller than the diameter of another cover $\beta$ ``qualifies" the relation $\alpha\to \beta$, i.e., $\beta$ is a finer subcover of $\alpha$. Corollary~\ref{cor:sdol}  implicitly relies on the following elementary result :

\begin{lemma} \label{lem:od30}
    Given two functors $F,G : X\to Y$ between preorders $X, Y$, there is a natural transformation $\eta : F\Rightarrow G$ iff $G$ is pointwise greater than $F$, i.e., for each $x\in ob(X)$, $F(x) \to G(x)$. In that case the natural transformation $\eta$ is unique.
\end{lemma}

Lemma~\ref{lem:od30} says that any existence result on natural transformation between functors is simply a statement about such pointwise inequalities.

\paragraph{A reinterpretation} Recall that post-composition with $F$ leads to a functor between functor categories
\[ F\circ : C^I \to D^I, \quad H \mapsto F\circ H. \]
Similarly, one has the functor $G\circ : C^I \to D^I$. These two functors $F\circ$ and $G\circ$ have the same domain and codomain, so one can construct their comma category  $\Comma{ F\circ }{ G\circ }$. Theorem~\ref{thm:qlf_actn} says that the qualifier $\Psi$ of the pair $(F,G)$ induces the map
\begin{equation} \label{eqn:qlf_actn}
    \Psi_* : \Comma{ F\circ }{ G\circ } \to \Hom_{C^I} 
\end{equation}
Equation~\eqref{eqn:qlf_actn} summarizes the primary role played by qualifying pairs. Note that the elements in $\Hom_{C^I}$ are natural transformations between functors from $I$ to $C$. This ability of qualifying pairs to produce natural transformations has important consequences in studying limiting behavior. We study this next.

\section{Qualifying pair and (co)limits} \label{sec:qpl}

Our final goal is a categorical understanding of entropy, and that involves interpreting many asymptotic dynamic properties as colimits. In this section we review how colimits change as functors are extended or adjointed.
We begin with a chain of functors $N \xrightarrow{P} X \xrightarrow{m} Z$ between some categories $N,Z$. Then in general one would have

\begin{lemma} \label{lem:df30}
For any sequence of functors $N \xrightarrow{P} X \xrightarrow{m} Z$, if the colimits $\colim m$ and $\colim mP$ exist, then there is a morphism $\colim mP \to \colim m$.
\end{lemma}

This has special significance if $Z$ is a preorder, in which one can interpret the relation $\colim mP \to \colim P$ as the colimit of $mP$ being greater than the colimit of $P$. If $m:X\to Z$ is interpreted as a measure, and $P:N\to X$ a parameterized family in $X$, then according to Lemma~\ref{lem:df30}, the colimit of the measure on the parameterized family is greater than the colimit of the measure. However, if the parameterization $P$ bears a certain spanning property, then the colimits would be equal. We now formally define this spanning via a generalization of inverses.

\paragraph{Pre- and post- right adjoints} Let $T:A\to B$ be a functor. Then $T^* : B\to A$ is said to be a \emph{ pre-right adjoint } to $T$ if $TT^* \Rightarrow \Id_B$. Such a functor $T$ is said to be \emph{pre-right adjoint enabled} or in brief \emph{pre-r.a.e.}. 
Similarly, $T^* : B\to A$ is said to be a \emph{ post-right adjoint } to $T$ if $\Id_B \Rightarrow TT^*$. In that case $T$ will be called \emph{post-right adjoint enabled} or in brief \emph{post-r.a.e.}. 

Post-right adjoints generalize the notion of divergence to infinity of sequences. As an example, take $A = \Nplus$ and $B=\Rplus$. Any functor $T:A\to B$ corresponds to a non-decreasing sequence of non-negative numbers $0 \leq x_1 \leq x_2 \leq \ldots$. Suppose 
\[ \bar{x} := \colim T = \lim_{n\to\infty} x_n = \infty .\]
Now define a map $T^* : \real^+_0 \to \num$ as 
\[ T^*(x) := \min \SetDef{ n\in \num }{ x_n \geq x } . \]
Since $T^*$ is order preserving it is a functor $T^* : \Rplus \to \Nplus$. Note that by construction, we have $TT^* \Rightarrow \Id_{\Rplus}$. Also, it was possible to define $T^*$ only because the $x_n$s converge to infinity. Thus although $T : \Nplus \to \Rplus$ is an embedding, it \emph{spans} the whole range of $\Rplus$. This idea of spanning is captured by a pre-right adjoint. The following lemma shows an important consequence of this spanning action.

\begin{lemma} \label{lem:sd40}
[Limit and colimit preservation] Let  $R:B\to C$ be a functor. 
\begin{enumerate} [(i)]
    \item Suppose $T:A\to B$ has a pre-right adjoint. Then if the limits $\lim R$ and  $\lim RT$ exist, then they are equal.
    \item Suppose $T:A\to B$ has a post- right adjoint. Then if the colimits $\colim R$ and $\colim RT$ exist, then they are equal.
\end{enumerate}
\end{lemma}

Lemma~\ref{lem:sd40} is proved in Section~\ref{sec:proof:sd40}. We are now ready to state our main theorem.

\begin{theorem}[Limits under qualifying pair] \label{thm:lim_qual}
    Let $Z$ be a co-complete category, and $(D,L)$ is a qualifying pair of functors arranged in the following diagram of categories and functors
    \[\begin{tikzcd} [column sep = huge]
        N \arrow{r}{P} & X \arrow{dd}[swap]{D} \arrow[dashed]{drrr}{L} \arrow{r}{m} & Z \\
        & & & & \mathcal{E} \\
        & \mathcal{E} \arrow{r}{ \REnv{m}{D} } & Z
    \end{tikzcd}\]
    Then :
    \begin{enumerate} [(i)]
        \item There is a morphism $\colim mP \to \colim \REnv{m}{D}$.
    \end{enumerate}
    Henceforth we assume a functor $P:N\to X$ such that $DP : N\to \mathcal{E}$ has a post-right adjoint.
    \begin{enumerate} [(i), resume] 
        \item the functor $P$ has a post-right adjoint.
        \item $\colim m = \colim mP$.
        \item $\colim \REnv{m}{D} \to \colim mP$. 
        \item If $Z$ is a preorder, then all the above morphisms are equalities.
    \end{enumerate}
\end{theorem}

Theorem~\ref{thm:lim_qual} is proved in Sec~\ref{sec:proof:lim_qual}. Theorem~\ref{thm:lim_qual} says that a parameterized limit of a measurement $m$ would lead to the same limit as $m$, provided the parameterization ``spans the end" of the space. This notion of spanning the end is captured in the requirement of having a post-right adjoint. 
Theorems \ref{thm:qualify} and \ref{thm:lim_qual} lead to this immediate result :

\begin{corollary}[Colimits along sequences] \label{cor:ddkw}
    Now suppose that $X$ has finite coproducts, $\mathcal{E}$ has countable coproducts, and $\paran{ x_n }_n$ is a sequence of objects in $X$ such that the sequence $D(x_n)$ is monotonic and has a post-right adjoint. Let $m:X\to Z$ be a functor into a co-complete category $Z$. Then $\colim \REnv{m}{D} \to \limsup_{n\to \infty} m(x_n)$.
\end{corollary}
Corollary~\ref{cor:ddkw} is proved in Section~\ref{sec:proof:ddkw}.
This completes the statement of all our results from general Category theory. We next investigate some functors that arise from the dynamical system.

\section{Functors induced by the dynamics} \label{sec:dynamics}

We are now in a position to describe the categorical interpretation of classical entropy, as summarized in Figure~\ref{fig:outline1}. We begin with the first definition, which is based on open covers. The category $\OpenCov$ has finite coproducts. The category theoretic coproduct  of two open covers corresponds to the operation of ``join", defined as
\[ a \vee b := \SetDef{ S_a \cap S_b }{ S_a\in a, \, S_b\in b } . \]
Since $\OpenCov$ is also a preorder, a product of a finite collection of elements is equivalent to a minimum. Next note that the dynamics behaves as an endofunctor on $\OpenCov$. In general, coproducts and endofunctors combine to produce a \emph{dynamics functor} :

\begin{lemma} \label{lem:prk4}
Let $\mathcal{C}$ be any category with finite coproducts, and $F:\mathcal{C}\to \mathcal{C}$ be an endofunctor. Then this induces a functor $\dyn_{F} : \Nplus \times \mathcal{C} \to \mathcal{C}$, whose action on objects is given by
\[ \dyn_{F}(n, c) := \vee_{t=0}^{n-1} F^t c \]
\end{lemma}

The proof of Lemma~\ref{lem:prk4} is direct and will be omitted. We call $\dyn_{F}$ the dynamics functor associated to the endofunctor $F$. Now note that given any open cover $a\in \OpenCov$, one gets a new open cover
\[ \map^{*} a := \SetDef{ \map^{-1} S_a }{ S_a\in a } . \]
The fact that the elements in the collection $\map^{*} a$ are open sets follows from the continuity of $\map$. Thus we get a functor 
\[ \map^{*} : \OpenCov \to \OpenCov \]
and as a result, by Lemma~\ref{lem:prk4}, the associated dynamics functor :
\begin{equation} \label{eqn:def:dyn}
    \dyn_{\map} : \Nplus \times \OpenCov \to \OpenCov, \quad (n,\alpha) \mapsto \alpha \wedge \map^{*} \alpha \wedge \cdots \wedge \map^{*(n-1)} \alpha.
\end{equation}
For brevity, we shall denote $\dyn_{\map}$ by $\dyn$. This is the point where we introduce the first notion of complexity. Consider the following map
\begin{equation} \label{eqn:cov}
    \CovNum : \OpenCov \to \Nplus, \quad C \mapsto \mbox{ minimum size of a subcover of } C.
\end{equation}
This map $\CovNum$ clearly preserves the partial order on $\OpenCov$ and is thus a functor too. Now consider the composite of functors 
\[\begin{tikzcd} \Nplus \times \OpenCov \arrow[bend left = 20, dashed]{rrrr}{ \Entropy } \arrow{rr}[swap]{ \dyn } && \OpenCov \arrow{rr}[swap]{\CovNum} && \Rplus \end{tikzcd}\]
This functor will be our first notion of complexity or \emph{entropy}. We will find it more useful to rewrite $\Entropy$ as
\begin{equation} \label{eqn:def:entropy}
    \Entropy := \CovNum \circ \dyn : \OpenCov \to \Rplus^{ \Nplus }
\end{equation}
This is due to the natural bijection between functors $\Nplus \times \OpenCov \to \Rplus$ and $\OpenCov \to \Rplus^{ \Nplus }$. The category $\Rplus^{ \Nplus }$ is of special importance to us, we call it the category of \emph{rates}. For brevity we denote it as
\begin{equation} \label{eqn:def:rates}
    \Z := \Rplus^{ \Nplus } .
\end{equation}
The category $\Z$ has as objects all non-decreasing sequences of non-negative numbers, ordered by pointwise $\leq$ relations. The functors $\Entropy$ and $\diam$ fit together along with Kan extensions into the following commuting diagram :
\[\begin{tikzcd} [column sep = large]
    \Z & \Z & \Z \\
    & \OpenCov \arrow[dashed]{ul}[name=x2]{} \arrow[dashed]{ur}[name=x1]{} \arrow{d}{\diam} \arrow{u}[name=F]{\Entropy} \\ 
    & \EpsCat \arrow[bend right=20]{uur}[swap]{ \REnv{\diam}{\Entropy} } \arrow[bend left=20]{uul}{ \LEnv{\diam}{\Entropy} } &
    \arrow[shorten <=2pt, shorten >=3pt, Rightarrow, to path={(F) -- (x1)},]{}
    \arrow[shorten <=2pt, shorten >=3pt, Rightarrow, to path={(x2) -- (F)},]{}
\end{tikzcd}\]
The following left Kan extension will help capture the first two notions of complexity in \eqref{eqn:pressure} :
\begin{equation} \label{eqn:def:P1_P4_fnctr}
    \Pressure_1 := \REnv{\diam}{ \Entropy } : \EpsCat \to \Z .
\end{equation}
This completes an task of representing the notions of covering complexity as a functor. 

\begin{figure}[!ht]\center	\includegraphics[width=.95\linewidth]{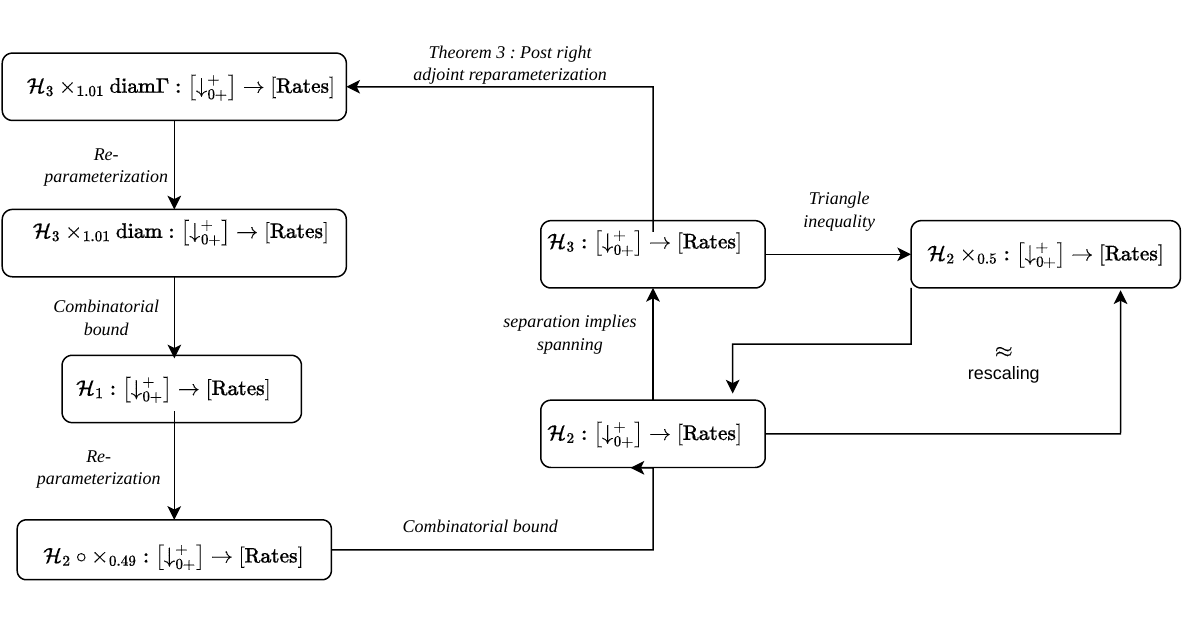}
	\caption{ Entanglement network between entropy functors. The figure shows a graph on 7 vertices, with each vertex drawn as a box. For each $1\leq i \leq 7$, the $i$-th vertex represents a functor $F_i : C_i \to \Z$, from some category $C_i$ into the same codomain $\Z$, which is the rate category \eqref{eqn:def:rates}. Each such functor represents a functorial / structural way of assigning a complexity sequence to a topological dynamical system $f:\Omega\to \Omega$. Each edge from a vertex $i$ to a vertex $j$ represents a sub-structural relationship between $F_i$ and $F_j$, made precise in Theorem~\ref{thm:entangle}. The existence of such an edge implies that the colimit of $F_i$ is less than or equal to the colimit of $F_j$. Since this entangle network is strongly connected, all these colimits must be equal. These colimits are related to the classical entropy limits expressed in \eqref{eqn:pressure}, and thus they are proved to be equal. }
    \label{fig:outline2}
\end{figure}

Next we provide functorial characterizations of $H_2, H_3$. It begins with the realization that the collection of finite subsets of $\Omega$ form a preorder $\FinSet$, with a preorder structure provided by the inclusion relations. Thus given two finite subset $A,A'$ of $X$, $A\to A'$ iff $A\subseteq A'$. As a result, cardinality is a functor :
\[ \cardinality : \FinSet \to \Rplus . \]
The quantities $\text{metric-sep}$ and $\text{metric-span}$ from \eqref{eqn:dpd9} can now be stated as functors
\begin{equation} \label{eqn:def:OrbCovSep}
    \begin{split}
        \OrbitSpan : \Nplus \times \FinSet \to \EpsCat, & \quad (n,A) \mapsto \text{metric-span} \paran{ A, d_n } \\
        \OrbitSep : \Nplus \times \FinSet \to \EpsCat, &\quad (n,A) \mapsto \text{metric-sep} (A, d_n) 
    \end{split}
\end{equation}
The functors $\OrbitSpan$ and $\OrbitSep$ are functors naturally induced by the dynamics, along with the functor $\Entropy$ from \eqref{eqn:def:entropy}. Given a pair $(n,A)$ in $\mathbb{N}\times \FinSet$, the notation $\tilde{ \norm{ (n,A) } }$ will be interpreted as the composite functor :
\[\begin{tikzcd} [column sep = huge] 
\Nplus \arrow[bend left = 20]{rr}{ \tilde{\cardinality} } \times \FinSet \arrow{r}[swap]{\proj_2} & \FinSet \arrow{r}[swap]{ \cardinality } & \Rplus
\end{tikzcd}\]
Next we interpret the operations of taking extrema along one variable in \eqref{eqn:def:P2_P3}, into Kan extensions
\[\begin{tikzcd} [column sep = huge]
    \Nplus \times \EpsCat \arrow{d}[swap]{ \Pressure_2 } & \Nplus \times \FinSet \arrow[dashed]{dl}[ name=x1 ]{} \arrow[dashed]{dr}[ name=x3 ]{} \arrow{d}[name=x2]{ \tilde{
    \cardinality} } \arrow{l}[swap]{ \Id_N \times \OrbitSpan } \arrow{r}{ \Id_N \times \OrbitSep }  & \Nplus \times \EpsCat \arrow{d}{ \Pressure_3 } \\
    \Rplus & \Rplus & \Rplus
    \arrow[shorten <=5pt, shorten >=30pt, Rightarrow, to path={(x1) to[out=-10,in=190] (x2)} ]{  }
    \arrow[shorten <=5pt, shorten >=30pt, Rightarrow, to path={(x2) to[out=-10,in=190] (x3)} ]{  }
\end{tikzcd}\]
These functors $\Pressure_2$ and $\Pressure_3$ reveal the connection of the analogous quantities from \eqref{eqn:def:P2_P3} with the topological structure of the dynamics. In summary :
\begin{equation} \label{eqn:def:_P2_P3_fnctr}
    \begin{split}
        \Pressure_2 := \LEnv{ \Id_N \times \OrbitSpan }{ \abs{\cdot} \circ \proj_2 }, &\quad \Pressure_2(n,\epsilon) := \inf \SetDef{ \abs{A} }{ A \in \FinSet, \, \OrbitSpan(A,n) > \epsilon } \\ 
        \Pressure_3 := \REnv{ \Id_N \times \OrbitSep }{ \abs{\cdot} \circ \proj_2 }, &\quad \Pressure_3(n,\epsilon) := \sup \SetDef{ \abs{A} }{ A \in \FinSet, \, \epsilon > \OrbitSep(A,n) } 
    \end{split}
\end{equation}
Finally, we recast the functors $\Pressure_2, \Pressure_3 : \Nplus \times \EpsCat \to \Rplus$ as 
\[ \Pressure_2 : \EpsCat \to \Z, \quad \Pressure_3 : \EpsCat \to \Z ,\]
where $\Z$ was defined in \eqref{eqn:def:rates}. This completes the task of expressing each of the notions of entropy as functors. In the next section we examine how they are inter-related via natural transformations. An overview of their inter-relations is presented in Figure~\ref{fig:outline2}.

\section{Network of diagrams} \label{sec:network}

In this section we prove Corollary~\ref{cor:thermo} using the category theoretical machinery we have been building. The following theorem is the key tool using which we interconnect the various functors developed so far.

\begin{theorem}[Entanglement theorem] \label{thm:entangle}
Suppose $S$ is a cocomplete category, and $\mathcal{G}$ is a directed graph on n vertices. Moreover
    \begin{enumerate}[(i)]
        \item For each vertex $i$, there is a category $C_i$ and a functor $F_i : C_i \to S$.
        \item For each edge $i\to j$, there is a category $A_{i,j}$, a functor $F_{i,j} : C_i \to C_j$ with a post right adjoint, a functor $G_{i,j} : A_{i,j} \to C_j$ and a natural transformation $\eta_{i,j}$ as shown
        \begin{equation} \label{eqn:pkd4}
            \begin{tikzcd} [column sep = large]
                C_j \arrow{d}[swap]{F_j} & A_{i,j} \arrow[dashed]{dl}[name=x2]{} \arrow{l}[swap]{ G_{i,j} } \arrow[dashed]{dr}[name=x1]{} \arrow{r}{F_{i,j}} & C_i \arrow{d}{F_i} \\
                S & & S 
                \arrow[shorten <=10pt, shorten >=1pt, Rightarrow, to path={(x1) to[out=225,in=-45] node[yshift=-10pt]{$\eta_{i,j}$} (x2)} ]{  }
        \end{tikzcd}
        \end{equation}
    \end{enumerate}
    Then for every $1\leq i,j \leq n$,
    \begin{enumerate} [(i)]
        \item If there is a path from i to j, then $\colim F_i \to \colim F_j$.
        \item In particular, if $S$ is a preorder, then $i,j$ are in the same connected component of $\mathcal{G}$ iff $\colim F_i = \colim F_j$.
    \end{enumerate}
\end{theorem}

Theorem~\ref{thm:entangle} is proved in Section~\ref{sec:proof:limit_ntwrk}. 

\paragraph{Remark} An analogous result can be stated for limits instead of colimits, in which the pre-right adjoint requirement has to replaced by an analogous pre-left adjoint condition. 

\paragraph{Remark} The natural transformations $\eta_{i,j}$ play the role of an \emph{entanglement}. This entanglement is not between functors $F_i$ and $F_j$, but between their extensions from a \emph{pre-parameterization} category $A_{i,j}$. 

\paragraph{Existence of post-right adjoints} When applying Theorem~\ref{thm:entangle} to any arrangement of functors and categories, one would need to verify the condition that the functor $F_{i,j}$ has a post-right adjoint. The following lemma provides a tool for verifying this.

\begin{lemma}[Composition of post-r.a.e. functors] \label{lem:ddp3}
    The composition of two functors with post right adjoints, also has a post right adjoint.
\end{lemma}

Lemma~\ref{lem:ddp3} is proved in Section~\ref{sec:proof:ddp3}. One key functor we shall rely on is defined the functorial version of the multiplication by $a$ :
\begin{equation} \label{eqn:def:times_a}
    \times_{a} : \EpsCat \to \EpsCat, \quad x \mapsto \begin{cases} ax & x \leq 1/a \\ 1 & x \geq 1/a \end{cases} .
\end{equation}
It is a standard practice in calculus to rely on such scaling functions to preserve the study of limits as $\epsilon\to 0^+$. This is formalized below.

\begin{lemma} \label{lem:lp4d9}
    For every positive number $a>0$, the functor $\times_{a} : \EpsCat \to \EpsCat$ has a post right adjoint, which is division by $a$. Further, division by $a$ is a right inverse to $\times_a$. 
\end{lemma}

We shall use a combination of Lemmas \ref{lem:ddp3} and \ref{lem:lp4d9} to interpret the relations between the various entropy functors in the light of Theorem \ref{thm:entangle}.

\paragraph{Structural relations} We list three natural transformations that are fundamental to the functors induced by the dynamics. These provide a precise meaning to the edge connections between the various complexity functors, as drawn in Figure~\ref{fig:outline2}.

\begin{enumerate}
    \item The following commuting diagram is a functorial restatement of metric relations \citep[][Sec 7.2, Rem 5]{Walters2000},
    \begin{equation} \label{eqn:Sec7_2_Rem5}
        \begin{tikzcd} [column sep = large]
             & \Z & \Z \\
            \Z & \EpsCat \arrow[dashed, bend right = 20]{ur}[name=x3]{} \arrow{l}[name=x1]{ \Pressure_2 } \arrow{u}[name=x2, swap]{ \Pressure_3 } \arrow{r}[swap]{ \times 0.5 } & \EpsCat \arrow{u}[swap]{ \Pressure_2 }
            \arrow[shorten <=10pt, shorten >=10pt, Rightarrow, to path={(x1) to[out=90, in=180] node[yshift=10pt]{$\eta_{2,3}$} (x2)} ]{  } 
            \arrow[shorten <=10pt, shorten >=2pt, Rightarrow, to path={(x2) to[out=0,in=135] node[yshift=10pt]{$\eta_{3,2}$} (x3)} ]{  }
        \end{tikzcd} 
    \end{equation}
    The transformation $\eta_{3,2}$ simply says that the counts provided by metric span is greater than the counts provided by metric separation. On the other hand, $\eta_{2,3}$ reverses this inequality, if one halves the grain size $\epsilon$.
    \item The following commuting diagram is  a functorial restatement of combinatorial relations     \citep[][Thm 7.7]{Walters2000} between covering complexity, and metric span and separation.
    \begin{equation} \label{eqn:Corr7_7_1}
        \begin{tikzcd} [column sep = large]
            \EpsCat \arrow{d}{\Pressure_3} & \EpsCat \arrow{l}[swap]{ \times 1.01 } & \OpenCov \arrow[dashed]{dll}[name=x1]{} \arrow{l}[swap]{\diam} \arrow[dashed]{drr}[name=x3]{} \arrow{d}[name=x2]{ H } \arrow{r}{\LebNum} & \EpsCat \arrow{r}{ \times 0.49 } & \EpsCat \arrow{d}{\Pressure_2} \\
            \Z & & \Z & & \Z
            \arrow[shorten <=10pt, shorten >=10pt, Rightarrow, to path={(x1) to[out=-10,in=180] node[yshift=-10pt]{$\eta_{3,H}$} (x2)} ]{ } 
            \arrow[shorten <=10pt, shorten >=5pt, Rightarrow, to path={(x2) to[out=-10,in=225] node[yshift=-10pt]{$\eta_{H,2}$} (x3)} ]{  } 
        \end{tikzcd}
    \end{equation}
    \item The final structural stems from purely topological properties. Now suppose $C$ is an UDC of diameter / grain $\epsilon$. Consider the UDC $\tilde{C}$ with disks having the same centers as those of $C$, but with radius $2.01\epsilon$. Now consider any point $x\in \Omega$. Then there must be a point $y$ such that the disk $B(y, \epsilon)$ of diameter $\epsilon$ and center $y$ lies in the cover $C$. Then note that by the triangle inequality, the disk $B(x, \epsilon)$ must lie within the disk $B(y, 2.1\epsilon)$. Therefore, the UDC $\tilde{C}$ has a Lebesgue number at least $\epsilon$ and is coarser than $C$. Thus we have a diagram : 
    \begin{equation} \label{eqn:kkd09}
        \begin{tikzcd} [column sep = large]
            & & \EpsCat \\
            \EpsCat \arrow[ bend left = 10 ]{urr}{ \times_{2.01} } & \udc \arrow[dashed]{dr}[name=x1]{} \arrow{l}[swap, name=x2]{ \Grain } \arrow{r}{\Expand} & \udc \arrow{u}[swap]{ \diam } \arrow{d}{\Leb} \\
            & & \EpsCat
            \arrow[shorten <=10pt, shorten >=3pt, Rightarrow, to path={(x1) to[out=180,in=-90] node[yshift=-10pt]{$\eta_{exp}$} (x2)} ]{  }
        \end{tikzcd}
    \end{equation}
    Here $\Expand$ is the functorial expansion from $C$ to $\tilde{C}$ just described, and $\times_{2.01} : \EpsCat \to \EpsCat$ is defined as in \eqref{eqn:def:times_a}. In fact since $\Grain = \diam \circ \iota_{UDC}$, we have from \eqref{eqn:kkd09} :
    \[ \Leb \Expand \Rightarrow \diam \iota_{UDC} , \]
    which by the qualifying action of $\paran{ \diam, \Leb }$ implies $\Expand \Rightarrow \Id_{UDC} $. 
\end{enumerate}


\paragraph{Exhaustive sequence of covers} In our application of Theorem~\ref{thm:entangle}, in most cases our choice of $A_{i,j}$ will be $\Nplus$. Note that we have at our disposal a functor $\Gamma : \Nplus \to \OpenCov$ such that in the following diagram
    \begin{equation} \label{eqn:Cdiam}
        \begin{tikzcd}
            \Nplus \arrow{r}{\Gamma} \arrow{drr}[swap]{ \diam \Gamma} & \udc \arrow{r}{ \iota_{UDC} } & \OpenCov \arrow{d}{\diam} \\ 
            & & \EpsCat
        \end{tikzcd}
    \end{equation}
    the functor $\Gamma \diam$ has a post-right adjoint. This is because one can find a nested sequence of UDC-s $\braces{ \Gamma_n }_{n=1}^{\infty}$ such that their diameters monotonically converge to $0^+$. Set $ \tilde{\Gamma}_n := \Expand \paran{ \Gamma_n } $, the UDC of grain $2 \diam \paran{ C_n }$.

 We now combine the various functors and natural transformations from \eqref{eqn:Sec7_2_Rem5}, \eqref{eqn:Corr7_7_1}, \eqref{eqn:kkd09}, \eqref{eqn:Cdiam} into the following \emph{entanglement} diagram :
\begin{equation} \label{eqn:entangle}
\begin{tikzcd} [column sep = huge, scale cd = 1.2]
\input{entangle.tex}
\end{tikzcd}
\end{equation}
The functors shown in blue will play the role of the $\Pressure_i$ in Theorem~\ref{thm:entangle}. The natural transformations shown in red correspond to the entanglements $\eta_{i,j}$. The equality of limits in \eqref{eqn:pressure} involve limits of the form $\lim_{n\to\infty} \frac{1}{n} \ln \paran{\cdot}$. This corresponds to the following functor
\[ \LogRate : \Z \to \realcat^{ \Nplus }, \quad \paran{a_n}_{n\in\num} \mapsto \paran{ \frac{1}{n} \ln a_n }_{n\in\num} .\]
Now the category $\realcat^{ \Nplus }$ has as objects all non-decreasing sequences of reals. The action of taking supremum of such sequences corresponds to the colimit functor. Concatenating, we get the following functor shown by a dashed arrow :
\[\begin{tikzcd} [column sep = large]
   \Z \arrow[dashed]{dr}[swap]{ \LogLim } \arrow{r}{\LogRate} & \realcat^{ \Nplus } \arrow{d}{ \colim } \\
   & \realcat
\end{tikzcd}\]
This functor $\LogLim$ has been appended to the entanglement diagram in \eqref{eqn:entangle}, and indicated in orange.

\paragraph{Proof of Corollary \ref{cor:thermo}} Consider an instance of Theorem~\ref{thm:entangle} in which the graph has three nodes. The functors $\Pressure_1, \Pressure_2, \Pressure_3$ have already been defined. See Table~\ref{tab:entangle} for a list of the functor $C_i A_{i,j} \xleftarrow{ F_{i,j} } \xrightarrow{ G_{i,j} } C_j $, and the natural transformations $\eta_{i,j}$.

\begin{table}
\caption{Entanglements between various entropy functors. The entanglement is in terms of Theorem~\ref{thm:entangle}. The various functors that form the components of this entanglement are given in the diagram in \eqref{eqn:entangle}.}
\begin{tabularx}{\linewidth}{|L|L|L|L|L|}
\hline
Edge $i\to j$ & $A_{i,j}$ & $F_{i,j}$ & $G_{i,j}$ & $\eta_{i,j}$  \\ \hline
$2\to 1$ & $\Nplus$ & $\times_{1.01} \diam \iota_{UDC} C$ & $\iota_{UDC} C$ & $\eta_{2,H}$ \\ \hline
$1\to 3$ & $\Nplus$ & $\iota_{UDC} \tilde{C}$ & $\times_{0.49} \LebNum \iota_{UDC} \tilde{C}$ & $\eta_{H,3}$ \\ \hline
$3\to 2$ & $\EpsCat$ & $\Id$ & $\Id$ & $\eta_{3,2}$ \\ \hline
\end{tabularx}
\label{tab:entangle}
\end{table}

It only remains to verify that the given $F_{i,j}$ are indeed with post-right adjoints. By construction, $\diam\circ \iota_{UDC} \circ C$ is post-r.a.e.. By Lemma~\ref{lem:ddp3} and Lemma~\ref{lem:lp4d9}, $F_{2,1}$ is post-r.a.e. too. By a similar reasoning, the functor $\times_{2.01} \circ \diam\circ \iota_{UDC} \circ C$ is post-r.a.e. too. But by commutation in \eqref{eqn:entangle}, this functor equals $\diam \circ \iota_{UDC} \circ \tilde{C}$ is post-r.a.e. too. Since $\paran{\diam, \LebNum}$ is a qualifying pair of functors by Theorem~\ref{thm:qualify}, by Theorem~\ref{thm:lim_qual}~(iii), $F_{1,3} = \iota_{UDC} \circ \tilde{C}$ is post-r.a.e. too. Finally, $F_{3,1}$ is trivially post-r.a.e.. Thus the nodes $1,2,3$ belong to the same strongly connected component of the entanglement graph, so by Theorem~\ref{thm:entangle},
\begin{equation} \label{eqn:final1}
    \colim \Pressure_2 = \colim \Pressure_3 = \colim \Pressure_1 = \colim \REnv{ \diam }{ \Entropy },
\end{equation}
with the last equality following from Theorem~\ref{thm:lim_qual}~(iv). Equation \eqref{eqn:final1} is the functorial statement of Corollary \ref{cor:thermo}, which we are now ready to prove.  As a result we have :
\[\begin{split}
    & \DarkGreen{ \sup_{\alpha} \lim_{n\to\infty} \frac{1}{n} \ln \CovNum \paran{ \alpha^{(n)} } } = \colim \LogLim \Entropy, \quad \mbox{ categorical interpretation } \\
    & = \colim \LogLim \REnv{ \diam }{ \Entropy } = \colim \LogLim \Entropy \Gamma, \quad \mbox{ by Theorem~\ref{thm:lim_qual} }, \\
    & = \DarkGreen{ \lim_{\epsilon\to 0^+} \sup_{ \diam(\alpha) \geq \epsilon } \frac{1}{n} \ln \CovNum \paran{ \alpha^{(n)} } }, \quad \mbox{ categorical interpretation } \\
    &= \colim \LogLim \Pressure_2 = \colim \LogLim \Pressure_3 , \quad \mbox{ Theorem~\ref{thm:entangle} }, \\
    &= \DarkGreen{ \lim_{\epsilon\to 0^+} \lim_{n\to\infty} \frac{1}{n} \ln H_2(\epsilon, n) } = \DarkGreen{ \lim_{\epsilon\to 0^+} \lim_{n\to\infty} \frac{1}{n} \ln H_3( \epsilon, n ) }.
\end{split}\]
All the limits from \eqref{eqn:pressure} have been shown above in green. We have utilized their categorical interpretations, namely --  $\Pressure_1$ from 
\eqref{eqn:def:P1_P4_fnctr}, and $\Pressure_2, \Pressure_3$ from \eqref{eqn:def:_P2_P3_fnctr}, to establish their equality via the equality of colimits in \eqref{eqn:final1}. This proves Corollary \ref{cor:thermo}. \qed

This completes the description of our results. The rest of the paper contains the proofs to the various category theoretical theorems and lemmas.

\section{Proof of the theorems and lemmas} \label{sec:proofs}

\subsection{Proof of Lemma~\ref{lem:Leb_diam}} \label{sec:proof:Leb_diam}

The topological definition of the diameter of an open cover $c$ is the maximum diameter of any of its components. The diameter of a set is defined to be the maximum possible distance between any two elements of the set. Thus if a set $S$ has diameter $r$, then for every $r'>r$, it would fit within a disk of diameter $r'$. Conversely, one could define the diameter of $S$ to be the infimum of all those $r'$ for which there is a disk of diameter $r'$ that would contain $S$. Thus one may define the diameter of an open cover $c$ as
\[ \diam(c) := \inf \SetDef{ r }{ \mbox{ there is an UDC of diameter } r \mbox{ coarser than } c } . \]
But the collection of UDC-s which are coarser than $c$ are precisely the left-slice category $\Comma{ \iota_{UDC} }{c}$. Because of the ordering of the numbers $(0,1]$ in $\EpsCat$, the operation of taking infimum of a set becomes a colimit. Thus the diameter becomes
\begin{equation} \label{eqn:def:diam}
    \diam(c) = \colim \SqBrack{ \Comma{ \iota_{UDC} }{c} \xrightarrow{ \Forget } \udc \xrightarrow{ \Grain } \EpsCat } .
\end{equation}
This is precisely the pointwise definition of the left Kan extension, so indeed $\diam = \REnv{ \Grain }{ \iota_{UDC} }$.

The topological definition of the Lebesgue number of an open cover $c$ is the supremum of all $\delta>0$ such that given any disk $U$ of diameter $\delta$, there is a piece of $c$ which properly contains $U$. This is equivalent to saying that the open cover $\Free(r)$ formed by taking the union of all possible disks of diameter $r$, is finer than $c$. Thus
\[ \LebNum(c) = \sup \SetDef{r}{ \Free(r) \mbox{ is finer than } c } \]
To say that $\Free(c)$ is finer than $c$ is to say that there is a morphism $c \to \iota(\Free(r))$. Thus the Lebesgue number can be rewritten as
\begin{equation} \label{eqn:def:Leb}
    \LebNum(c) := \lim \SqBrack{ \Comma{ c }{ \iota \Free } \xrightarrow{\Forget} \EpsCat } .
\end{equation}
This is precisely the pointwise definition of the right Kan extension. So $\LebNum = \LEnv{ \Id_{\mathcal{E}} }{\iota \Free}$. This completes the poof of Lemma~\ref{lem:Leb_diam}. \qed

\subsection{Proof of Theorem~\ref{thm:qualify}} \label{sec:proof:qualify}

We begin with part~(i) by showing that :
\begin{equation} \label{eqn:sjnd}
    D = \REnv{G}{\iota} = \REnv{ \Id_{\mathcal{E}} }{\iota G^r} .
\end{equation}
To see why, first note that $D$ is a left extension of $\Id_{\mathcal{E}}$ along $\iota G^r$, namely :
\[ D \iota G^r \Leftarrow G G^r = \Id_{\mathcal{E}} . \]
To show that $D$is a left Kan extension, it has to be shown that $D$ is minimal. So let $F$ be any other left extension of $\Id_{\mathcal{E}}$ along $\iota G^r$, i.e., $\Id_{\mathcal{E}} \Rightarrow F \iota G^r$. It has to be shown that $D\Rightarrow F$. Then
\[ F\iota = F\iota \Id_{ \tilde{X} } \Leftarrow F\iota F^r G \Leftarrow \Id_{ \mathcal{E} } G = G . \]
Thus $F$ is a left extension of $G$ along $\iota$. But since $D = \REnv{G}{\iota}$, $G\Rightarrow F$ as claimed.
The following lemma is about such left and right Kan extensions of the same functor pair.

\begin{lemma} \label{lem:d30l}
Consider a pair of functors $A \xleftarrow{\alpha} B \xrightarrow{\beta} C$ such that 
\begin{enumerate} [(i)]
    \item $C$ is a preorder;
    \item there is a full-functor $\gamma : C\to A$ such that $\gamma \circ \beta = \alpha$;
    \item the right and left Kan extensions $\LEnv{\beta}{\alpha}, \REnv{\beta}{\alpha}$ exist.
\end{enumerate}
Then $\REnv{\beta}{\alpha} \Rightarrow \LEnv{\beta}{\alpha}$. 
\end{lemma}

Lemma~\ref{lem:d30l} follows from arguments in coend calculus and the author is referred to \citep[][Sec X]{Maclane2013}. Lemma~\ref{lem:d30l} is satisfied in our case by the substitutions :
\[ \begin{tikzcd} [column sep = large] B \arrow{d}[swap]{\alpha} \arrow{r}{\beta} & C \arrow{dl}{\gamma} \\ A \end{tikzcd} \quad \longrightarrow \quad
\begin{tikzcd} [column sep = large] \mathcal{E} \arrow{d}[swap]{\iota G^r} \arrow{r}{ \Id_{\mathcal{E}} } & \mathcal{E}  \arrow{dl}{\iota G^r} \\ X \end{tikzcd} \quad \longrightarrow \quad
\begin{tikzcd} [column sep = large] \EpsCat \arrow{d}[swap]{ \iota_{UDC} \Free } \arrow{r}{ \Id } & \EpsCat \arrow{dl}{ \iota_{UDC} \Free } \\ \OpenCov \end{tikzcd} \]
Thus by Lemma~\ref{lem:d30l}, $D\Rightarrow L$, satisfying the first requirement to being a qualifying pair.

Now suppose that $\mathcal{E}$ is totally ordered. Any object in the comma category $\Comma{L}{D}$ is pair $x, y\in ob(X)$ such that $L(y) \to D(x)$. By construction
\[ D(x) = \inf E/x, \, E/x := \SetDef{ e\in \mathcal{E} }{ \iota G^r(e) \to x }, \quad L(y) := \sup y/E, \, y/E := \SetDef{ e\in \mathcal{E} }{ y \to \iota G^r(e) } \]
There is an increasing sequence of points $l_1 \leftarrow l_2 \leftarrow l_3 \leftarrow \ldots$ in $y/E$ converging to $L(y)$ from the left. Similarly, there is a sequence of points $d_1 \rightarrow d_2 \rightarrow d_3 \rightarrow \ldots$ in $E/x$ converging to $D(x)$ from the right. Thus for some $n\in \num$,
\[ D(x) \leftarrow d_n \leftarrow l_n \leftarrow L(y), \quad d_n \in E/x, \, l_n \in y/E . \]
Then by definition of the sets $E/x$, $y/E$,
\[ x \leftarrow \iota G^r (d_n) \leftarrow \iota G^r (l_n) \leftarrow y. \]
Thus there is a mapping $\Psi$ of $\Comma{L}{D}$ into $\Comma{ \Id_X }{ \Id_X }$. This makes $(L,D)$ a qualifying pair of functors. This completes the proof of Theorem~\ref{thm:qualify}. \qed

\subsection{Proof of Theorem~\ref{thm:qlf_actn} } \label{sec:proof:qlf_actn}

Consider any morphism $\paran{ x,y,\phi } \xrightarrow{ (f,g) } \paran{ x',y',\phi' }$ in $\Comma{F}{G}$. Then we have the following commutations and maps :
\begin{equation} \label{eqn:ofpe}
    \begin{tikzcd} (x,y) \arrow{d}{ (f,g) } \\ (x', y') \end{tikzcd} \;\xleftarrow{ \Forget }\;
    \begin{tikzcd}
    Fx \arrow{d}{Ff} \arrow{r}{ \phi } & Gy \arrow{d}{Gg} \\
    Fx' \arrow{r}{ \phi' } & Gy'
    \end{tikzcd} \;\xrightarrow{\Psi}\; \begin{tikzcd}
    x \arrow{d}{f} \arrow{r}{ \Psi(\phi) } & y \arrow{d}{g} \\
    x' \arrow{r}{ \Psi(\phi') } & y'
    \end{tikzcd}
    \;\xrightarrow{ \Forget }\; \begin{tikzcd} (x,y) \arrow{d}{ (f,g) } \\ (x', y') \end{tikzcd}
\end{equation}
Now suppose $\eta : FF' \Rightarrow GG'$. This means :
\[ \forall \begin{tikzcd} i \arrow{d}{ \gamma } \\ i' \end{tikzcd} : \quad 
\begin{tikzcd}
    FF'i \arrow{d}{FF'\gamma} \arrow{r}{ \eta_i } & GG'i \arrow{d}{GG'\gamma} \\
    FF'i' \arrow{r}{ \eta_{i'} } & GG'i'
\end{tikzcd}.\]
This means there is a morphism 
\[\begin{tikzcd} \paran{ F'i, G'i, \eta_i } \arrow{rr}{ \paran{F'\gamma, G'\gamma} } & & \paran{ F'i', G'i', \eta_{i'} } \end{tikzcd}\]
so by \eqref{eqn:ofpe}, we have the commutation
\[\begin{tikzcd} [column sep = large]
    F'i \arrow{d}[swap]{ F'\gamma } \arrow{r}{ \Psi \paran{ \eta_i } } & G'i \arrow{d}{G'\gamma} \\
    F'i' \arrow{r}{ \Psi \paran{ \eta_{i'} } } & G'i'
\end{tikzcd}\]
Thus there is indeed a natural transformation from $F'$ to $G'$, with $\Psi(\eta_i)$ being the component of this transformation at any point $i\in ob(I)$. This completes the proof of Theorem~\ref{thm:qlf_actn}. \qed

\subsection{Proof of Lemma~\ref{lem:sd40}} \label{sec:proof:sd40}

We only prove Claim (i), as the proof for Claim (ii) is analogous.

Consider the functors $R:B\to C, RT:A\to C, RTT^* : B\to C$, mapping into $C$. Let their limits be respectively $c_R, c_{RT}, c_{RTT^*}$. Thus there is a natural transformation $\eta^{(R)}$ such that $\paran{ c_R, \eta^{(R)} }$ is the minimal cone over the functor $R$. Thus $\eta^{(R)}$ is a natural transformation from the constant $c_R$-valued functor $B\to C$, to $R$. One can define $\eta^{(RT)}, \eta^{(RTT^*)}$ analogously. In summary
\[ \eta^{(R)} : c_{R} \Rightarrow R, \quad \eta^{(RT)} : c_{RT} \Rightarrow RT, \quad \eta^{(RTT^*)} : c_{RTT^*} \Rightarrow RTT^* . \]
Any cone over $R$ can be restricted to a cone over $RT$. Thus there is a unique morphism $\phi : c_R \to c_{RT}$ which pulls back the cone $\paran{ c_{RT}, \eta^{(RT)} }$ into the restriction of the cone $\paran{ c_R, \eta^{(R)} }$ to $RT$, namely : 
\begin{equation} \label{eqn:ef440}
    \forall \begin{tikzcd} a \arrow{d}[swap]{f} \\ a' \end{tikzcd} \in A \quad \mbox{ implies } \quad  
\begin{tikzcd} Ta \arrow{d}[swap]{Tf} \\ Ta' \end{tikzcd} \quad \mbox{ implies } \quad 
\begin{tikzcd} [column sep = large]
    RTa \arrow{dd}[swap]{RTf} &  \\ 
    & c_{RT} \arrow[bend left=20]{ul}[swap]{ \eta^{(RT)}_{a} } \arrow[bend right=20]{dl}{ \eta^{(RT)}_{a'} } & c_{R} \arrow{l}{\phi} \arrow[bend right=20]{ull}[swap]{ \eta^{(R)}_{Ta} } \arrow[bend left=20]{dll}{ \eta^{(R)}_{Ta'} } \\
    RTa' 
\end{tikzcd}
\end{equation}
Since $T^*$ is a pre-right adjoint of $T$, there is a natural transformation $\eta : TT^* \Rightarrow \Id_B$. Thus 
\begin{equation} \label{eqn:frr2}
    \forall \begin{tikzcd} b \arrow{d}[swap]{h} \\ b' \end{tikzcd} \in B \quad \mbox{ implies } \quad  
\begin{tikzcd} [column sep = large]
    R b \arrow{d}{Rh} & RTT^* b \arrow{d}{RTT^* h} \arrow{l}{ R \eta_b } \\
    R b' & RTT^* b' \arrow{l}{ R \eta_{b'} }
\end{tikzcd}
\end{equation}
We now prove that $c_{R}$ is also the limit for $RTT^*$. To see why, consider any cone $c_{RTT^*}$ over $RTT^*$, with connecting morphisms $\eta^{(RTT^*)}$. By \eqref{eqn:frr2}, this extend to a cone over the functor $R$ and we must have a unique morphism $\gamma : c_{RTT^*} \to c_{R}$ such that the following diagram commutes.
\[\forall \begin{tikzcd} b \arrow{d}[swap]{g} \\ b' \end{tikzcd} \in B , \quad
\begin{tikzcd} [column sep = large]
    & Rb \arrow{dd}{Rg} & RTT^*b \arrow{l}[swap]{ R\eta_{b} }  \arrow{dd}[swap]{RTT^*g} &  \\ 
    c_{R} \arrow[bend left=20]{ur}{ \eta^{(R)}_{b} } \arrow[bend right=20]{dr}[swap]{ \eta^{(R)}_{b'} } & & & c_{RTT^*} \arrow[bend right=120]{lll}[swap]{\gamma} \arrow[bend left=20]{ul}[swap]{ \eta^{(RTT^*)}_{b} } \arrow[bend right=20]{dl}{ \eta^{(RTT^*)}_{b'} } \\
    & Rb' & RTT^*b' \arrow{l}{ R\eta_{b'} } 
\end{tikzcd}\]
The limiting cone with apex $c_{RT}$ along with the morphism $\phi$ from \eqref{eqn:ef440} to get :
\begin{equation} \label{eqn:gfy4}
    \forall \begin{tikzcd} b \arrow{d}[swap]{g} \\ b' \end{tikzcd} \in B , \quad
\begin{tikzcd} [column sep = large]
    & Rb \arrow{dd}{Rg} & RTT^*b \arrow{l}[swap]{ R\eta_{b} }  \arrow{dd}[swap]{RTT^*g} &  \\ 
    c_{R} \arrow[out=-90, in=-90]{rrrr}{ \phi } \arrow[bend left=20]{ur}{ \eta^{(R)}_{b} } \arrow[bend right=20]{dr}[swap]{ \eta^{(R)}_{b'} } & & & c_{RTT^*} \arrow[bend right=120]{lll}[swap]{\gamma} \arrow[bend left=20]{ul}[swap]{ \eta^{(RTT^*)}_{b} } \arrow[bend right=20]{dl}{ \eta^{(RTT^*)}_{b'} } & c_{RT} \arrow{l}{\psi} \arrow[bend right=20]{ull}[swap]{ \eta^{(RT)}_{T^*b} } \arrow[bend left=20]{dll}{ \eta^{(RT)}_{T^*b'} } \\
    & Rb' & RTT^*b' \arrow{l}{ R\eta_{b'} } 
\end{tikzcd}
\end{equation}
The diagram in \eqref{eqn:gfy4} reveals that $c_R$ itself is the apex of a cone over $RTT^*$, with connecting morphisms $\eta_b^{(RTT^*)} \psi \phi$. Moreover, the arbitary cone at $c_{RTT^*}$ factorizes through this cone. Thus $c_R$ must be the limit of $RTT^*$. Since the cone over $RTT^*$ from $c_R$ itself factorizes via the cone from $c_{RT}$, we must have $c_R = c_{RT} = c_{RTT^*}$, as claimed. \qed

\subsection{Proof of Theorem~\ref{thm:lim_qual}} \label{sec:proof:lim_qual}

Claim~(i) follows from the following chain of morphisms in $Z$ : 
\begin{equation} \label{eqn:hhp3}
    \begin{tikzcd} [column sep = large]
    \colim mP \arrow[dashed]{d}{} \arrow{rr}{ \mbox{ by Lem \ref{lem:df30} } } & & \colim m \arrow{d}{ \mbox{by \eqref{eqn:L_R_Env}} } \\
    \colim \REnv{m}{D} && \colim \paran{ \REnv{m}{D} \circ D } \arrow{ll}{ \mbox{ by  Lem \ref{lem:df30} } } .
\end{tikzcd}
\end{equation}
Next Claim~(ii) will be proven. Let $\Psi$ be the post-right adjoint of $DP$, i.e., $\Id_N \Rightarrow DP\Psi$. Thus
\[ L \Rightarrow \REnv{L}{D} D \Rightarrow DP\Psi \REnv{L}{D} D \]
$(L,D)$ form a qualifying pair of functors, so by Theorem~\ref{thm:qlf_actn}
\begin{equation} \label{eqn:sj3}
    \Id_X \Rightarrow P\Psi \REnv{L}{D} D
\end{equation}
Equation \eqref{eqn:sj3} provides the post-right adjoint of $P$ namely $\Psi \REnv{L}{D} D$. This proves Claim~(ii). Claim~(iii) now follows from Lemma~\ref{lem:sd40}~(ii).

To prove Claim~(iv) we compose both sides of \eqref{eqn:sj3} with the functor $m$ to get :
\[ m \Rightarrow m P \Psi \REnv{L}{D} D \]
By the minimality of the relation $m \Rightarrow \REnv{m}{D} D$, we must have
\begin{equation} \label{eqn:sd3jk}
    \REnv{m}{D} \Rightarrow m P \Psi \REnv{L}{D}
\end{equation}
In consequence, we get the ordering :
\begin{equation} \label{eqn:kdb}
    \begin{tikzcd} [column sep = large]
        \colim \REnv{m}{D} \arrow{r}{ \mbox{ by \eqref{eqn:sd3jk} } } & \colim m P \Psi \REnv{L}{D} \arrow{r}{ \mbox{Lem \ref{lem:df30} } } & \colim mP.
    \end{tikzcd}
\end{equation}
Claim(v) follows from \eqref{eqn:hhp3}, \eqref{eqn:kdb}, and the fact that $Z$ is a preorder. \qed

\subsection{Proof of Corollary~\ref{cor:ddkw}} \label{sec:proof:ddkw}

Note that the sequence $x_n$ can be interpreted as the image of a functor $Q:\braces{\num} \to X$. Consider the left Kan extension of $Q$ along the inclusion $\iota_+ : \braces{\num} \to \Nplus$, given by
\[ \bar{Q} := \REnv{ \iota_+ }{ Q } : \Nplus \to X, \quad \bar{Q}(n) := \coprod_{j=1}^{n} x_n .  \]
Also by assumption, there is a functor $R : \Nplus \to \mathcal{E}$ such that the following commutation holds :
\[\begin{tikzcd}
    \braces{\num} \arrow{d}[swap]{ \iota_+ } \arrow{r}{ Q } & X \arrow{d}{D} \\
    \Nplus \arrow{r}{ R } & \mathcal{E}
\end{tikzcd}\]
Therefore, we have 
\[ R \iota_+ = DQ \Rightarrow D \bar{Q} \iota_+ . \]
Since $\iota_+ : \braces{\num} \to \Nplus$ is bijective on objects, a natural transformation $R \iota_+ \Rightarrow D \bar{Q} \iota_+$ implies a natural transformation $R \Rightarrow D \bar{Q}$. Then since $R$ has a post right adjoint, say $\Psi : \mathcal{E} \to \Nplus$, we have
\[ \Id_{\mathcal{E}} \Rightarrow R \Psi \Rightarrow D \bar{Q} \Psi . \]
Thus $\Psi$ is a post right adjoint to $D\bar{Q}$ too. Thus, the functor $\bar{Q} : \Nplus \to X$ satisfies the criterion on the functor $P:N\to X$ from Claim (i). Thus we must have
\[ \colim \REnv{m}{D} \to \colim m \bar{Q} = \limsup_{n\to\infty} m(x_n), \]
proving Corollary~\ref{cor:ddkw}.

\subsection{Proof of Theorem~\ref{thm:entangle}} \label{sec:proof:limit_ntwrk}

Both claims of the theorem will be proved if it can be shown that for each edge $i\to j$ in $\mathcal{G}$, $\lim F_i \to \lim F_j$. Since $S$ is cocomplete, for any category $A$, the operation of taking colimits is a functor
\[ \lim : S^A \to S \]
Now take $A = C_i$. Then $F_i$ and $F_j \circ F_{i,j}$ are two objects in the functor category $S^{C_i}$, with $\eta_{i,j}$ being a morphism between them. Thus by the functorial property of $\colim$, we have $\colim F_i \to \colim F_j \circ F_{i,j}$. By Lemma~\ref{lem:df30}, $\colim F_j \circ F_{i,j} \to \colim F_j$ and thus $\colim F_i \to \colim F_j$. \qed

\subsection{Proof of Lemma \ref{lem:ddp3} } \label{sec:proof:ddp3}

Suppose that in the composite diagram $A \xrightarrow{F} B \xrightarrow{G} C$, both $F$ and $G$ have post right adjoints $F^*$ and $G^*$ respectively. Then note that
\[ \Id_C \Rightarrow GG^* = G \Id_B G^* \Rightarrow GFF^* G^* ,\]
making $F^*G^*$ a post right adjoint of $GF$. \qed

\section{Declarations}

\paragraph{Author’s Contribution} The sole author Suddhasattwa Das was responsible for all the research and in the preparation of the manuscript.
 
\paragraph{Conflict of Interest} There are financial or non-financial interests in competition with our work.

\paragraph{Availability of Data and Materials} The research is not dependent on any data, experimental or otherwise. Please contact the author Suddhasattwa Das for any further information.

\paragraph{Funding} No funding was obtained for this study.

\bibliographystyle{unsrt}
\bibliography{\Path References,Ref}
\end{document}

%% file: Paper_keywords.tex
\DeclareMathOperator{\Pressure}{\mathcal{P}}

\DeclareMathOperator{\Grain}{ Grain }
\DeclareMathOperator{\LebNum}{ Leb }
\DeclareMathOperator{\dyn}{ Dyn }
\DeclareMathOperator{\CovNum}{ cov }
\DeclareMathOperator{\OrbitSep}{ sep }
\DeclareMathOperator{\OrbitSpan}{ span }
\newcommand{\Z}{\mathcal{Z}}
\newcommand{\map}{ \mathcal{T} }
\newcommand{\Entropy}{h}
\DeclareMathOperator{\Expand}{Expnd}
\renewcommand{\REnv}[2]{ \mathcal{L}_{#2} \left( #1 \right) }
\renewcommand{\LEnv}[2]{ \mathcal{R}_{#2} \left( #1 \right) }
\newcommand{\LogLim}{\text{ LogLim }}

%% file: entangle.tex
    & & \realcat & \realcat & \realcat \\
& & \Z \arrow[orange]{u}{ \LogLim } & \Z \arrow[orange]{u}{ \LogLim } & \Z \arrow[orange]{u}[swap]{ \LogLim } \\
\Nplus \arrow[green, dashed, bend right=10]{drr}{\eta} \arrow{r}{ \Gamma } \arrow[dashed, bend right=20]{rdd}[swap]{ \tilde{\Gamma} } & \udc \arrow[bend right=10]{dd}{ \hspace*{-0.7cm}\rotatebox{90}{Expand} } \arrow{r}{ \iota_{UDC} } & \OpenCov \arrow{r}[name=D]{ \diam } \arrow[blue]{u}[name=P1a]{ \Pressure_1 := \Entropy } & \EpsCat \arrow{r}[swap]{ \times 1.01 } \arrow[dashed]{u}[name=P2]{} \arrow[bend left = 15]{dl}[swap]{ \times 2.01 } & \EpsCat \arrow[blue]{ul}[name=P2b]{ \Pressure_2 } \arrow[blue]{u}[name=P3b, swap]{ \Pressure_3 } \\
& & \EpsCat & & \\
& \udc \arrow{r}{ \iota_{UDC} } & \OpenCov \arrow{u}{ \diam } \arrow{r}[name=L]{ \LebNum } \arrow[blue]{d}[name=P1b, swap]{ \Pressure_1 := \Entropy } & \EpsCat \arrow{rd}{ \times 0.49 } \arrow[dashed]{d}[name=P3]{} & \\
& & \Z \arrow[orange]{d}[swap]{\LogLim} & \Z \arrow[orange]{d}[swap]{\LogLim} & \EpsCat \arrow[blue]{l}{ \Pressure_3 } \\
& & \realcat & \realcat & & 
\arrow[shorten <=10pt, shorten >=10pt, Rightarrow, to path={(P2) to[out=180, in=0] node[yshift=10pt]{$\eta_{2,H}$} (P1a)}, red ]{  }
\arrow[shorten <=10pt, shorten >=10pt, Rightarrow, to path={(P1b) to[out=0, in=180] node[yshift=-10pt]{$\eta_{H,3}$} (P3)}, red ]{  }
\arrow[shorten <=10pt, shorten >=10pt, Rightarrow, to path={(L) to[out=45, in=-45] node[xshift=15pt]{$\eta_{exp}$} (D)}, red ]{  }
\arrow[shorten <=1pt, shorten >=1pt, Rightarrow, to path={(P3b) to[out=180, in=45] node[yshift=10pt]{$\eta_{3,2}$} (P2b)}, red ]{  }

%% file: 2024_01_18.bbl
\begin{thebibliography}{10}

\bibitem{Walters1975}
P.~Walters.
\newblock A variational principle for the pressure of continuous
  transformations.
\newblock {\em Amer. J. Math.}, 97(4):937--971, 1975.

\bibitem{Walters2000}
P.~Walters.
\newblock {\em An introduction to ergodic theory}, volume~79.
\newblock Springer-Verlag New York, 2000.

\bibitem{PesinPitskel1984}
Y.~Pesin and B.~Pitskel.
\newblock Topological pressure and the variational principle for noncompact
  sets.
\newblock {\em Funct. Anal. Appl.}, 18(4):307--318, 1984.

\bibitem{Mummert2007}
A.~Mummert.
\newblock A variational principle for discontinuous potentials.
\newblock {\em Erg. Th. Dyn. Sys.}, 27(2):583--594, 2007.

\bibitem{Climenhaga_entropy}
V.~Climenhaga.
\newblock Entropy of {S}-gap shifts, 2014.

\bibitem{ClimenhagaThompson2012}
V.~Climenhaga and Daniel~J Thompson.
\newblock Intrinsic ergodicity beyond specification: $\beta$-shifts, s-gap
  shifts, and their factors.
\newblock {\em Israel J. Math.}, 192(2):785--817, 2012.

\bibitem{grassberger1991inf}
P.~Grassberger.
\newblock Information and complexity measures in dynamical systems.
\newblock In {\em Information dynamics}, pages 15--33. Springer, 1991.

\bibitem{vallee1987inf}
R.~Vall{\'e}e.
\newblock Information entropy and state observation of a dynamical system.
\newblock In {\em International conference on information processing and
  management of uncertainty in knowledge-based systems}, pages 403--405.
  Springer, 1987.

\bibitem{Spandl_sofic_2008}
C.~Spandl.
\newblock Computability of topological pressure for sofic shifts with
  applications in statistical physics.
\newblock {\em J. UCS}, 14(6):876--895, 2008.

\bibitem{BDWY2020}
M.~Burr, S.~Das, C.~Wolf, and Y.~Yang.
\newblock Computability of topological pressure on compact shift spaces beyond
  finite type.
\newblock {\em Nonlinearity}, 45:4250, 2022.

\bibitem{LongFerguson2019}
A.~Long and A.~Ferguson.
\newblock Landmark diffusion maps (l-dmaps): Accelerated manifold learning
  out-of-sample extension.
\newblock {\em Appl. Comput. Harmon. Anal.}, 47(1):190--211, 2019.

\bibitem{LiangPaisley2015}
D.~Liang and J.~Paisley.
\newblock Landmarking manifolds with {G}aussian processes.
\newblock In {\em International Conference on Machine Learning}, pages
  466--474. PMLR, 2015.

\bibitem{SilvaEtAl_landmark_2006}
J.~Silva, J.~Marques, and J.~Lemos.
\newblock Selecting landmark points for sparse manifold learning.
\newblock In {\em Advances in neural information processing systems}, pages
  1241--1248, 2006.

\bibitem{calcines2013limit}
J.~Calcines, L.~Paricio, and M.~Rodr{\'\i}guez.
\newblock Limit and end functors of dynamical systems via exterior spaces.
\newblock {\em Bull. Belgian Math. Soc.}, 20(5):937--959, 2013.

\bibitem{spivak2015steady}
D.~Spivak.
\newblock The steady states of coupled dynamical systems compose according to
  matrix arithmetic, 2015.

\bibitem{ngotiaoco2017}
T.~Ngotiaoco.
\newblock Compositionality of the {R}unge-{K}utta method, 2017.

\bibitem{jaz2020double}
D.~Myers.
\newblock Double categories of open dynamical systems, 2020.

\bibitem{Baez2014_bayesian}
J.~Baez and T.~Fritz.
\newblock A {B}ayesian characterization of relative entropy, 2014.

\bibitem{BaezEtAl2011_info}
J.~Baez, T.~Fritz, and T.~Leinster.
\newblock A characterization of entropy in terms of information loss.
\newblock {\em Entropy}, 13(11):1945--1957, 2011.

\bibitem{lomadze1999time}
V.~Lomadze.
\newblock Time-varying linear dynamical systems.
\newblock In {\em Proc. A. Razmadze Math. Inst}, volume 119, pages 121--132,
  1999.

\bibitem{Delvenne2019_dyn}
JC. Delvenne.
\newblock Category theory for autonomous and networked dynamical systems.
\newblock {\em Entropy}, 21(3):302, 2019.

\bibitem{Suda2022Poincare}
T.~Suda.
\newblock A categorical view of poincar{\'e} maps and suspension flows.
\newblock {\em Dynamical Systems}, 37(1):159--179, 2022.

\bibitem{MossPerrone2022ergdc}
S.~Moss and P.~Perrone.
\newblock A category-theoretic proof of the ergodic decomposition theorem.
\newblock {\em Erg. Th. Dyn. Sys.}, pages 1--27, 2022.

\bibitem{perrone2022kan}
P.~Perrone and W.~Tholen.
\newblock Kan extensions are partial colimits.
\newblock {\em Applied Cat. Struct.}, pages 1--69, 2022.

\bibitem{street2004categorical}
R.~Street.
\newblock Categorical and combinatorial aspects of descent theory.
\newblock {\em Applied Cat. Struct.}, 12(5):537--576, 2004.

\bibitem{Riehl_homotopy_2014}
E.~Riehl.
\newblock {\em Categorical homotopy theory}, volume~24.
\newblock Cambridge University Press, 2014.

\bibitem{polonsky2020local}
A.~Polonsky and P.~Johann.
\newblock Local presentability of certain comma categories.
\newblock {\em Applied Cat. Struct.}, 28(1):135--142, 2020.

\bibitem{winter2009arrow}
M.~Winter.
\newblock Arrow categories.
\newblock {\em Fuzzy Sets and Systems}, 160(20):2893--2909, 2009.

\bibitem{steingartner2014categorical}
W.~Steingartner and D.~Radakovi{\'c}.
\newblock Categorical structures as expressing tool for differential calculus.
\newblock {\em Open Computer Science}, 4(3):96--106, 2014.

\bibitem{grandis1997categorically}
M.~Grandis.
\newblock Categorically algebraic foundations for homotopical algebra.
\newblock {\em Applied Cat. Struct.}, 5(4):363--413, 1997.

\bibitem{burke2018synthetic}
M.~Burke.
\newblock A synthetic version of lie’s second theorem.
\newblock {\em Applied Cat. Struct.}, 26(4):767--798, 2018.

\bibitem{Maclane2013}
S.~Mac Lane.
\newblock {\em Categories for the working mathematician}, volume~5.
\newblock Springer Science \& Business Media, 2013.

\end{thebibliography}
